\documentclass[a4paper,12pt]{amsart}

\usepackage{amssymb}
\usepackage{xy}
\usepackage{bm}
\xyoption{all}
\xyoption{poly}

\newtheorem{theo}{Theorem}[section]
\newtheorem{coro}[theo]{Corollary}
\newtheorem{prop}[theo]{Proposition}
\newtheorem{lemma}[theo]{Lemma}
\newtheorem{conj}[theo]{Conjecture}
\theoremstyle{definition}
\newtheorem{defi}[theo]{Definition}
\newtheorem{remark}[theo]{Remark}

\def\cq{{\backslash \hskip -2pt \backslash}}

\def\CC{{\mathcal{C}}}
\def\CS{{\mathcal{S}}}

\def\CL{{\mathcal{L}}}

\def\CH{{\mathcal{H}}}

\def\CR{{\mathcal{R}}}

\def\CU{{\mathcal{U}}}
\def\CA{{\mathcal{A}}}

\def\CT{{\mathcal{T}}}
\def\CR{{\mathcal{R}}}
\def\CL{{\mathcal{L}}}

\def\CS{{\mathcal{S}}}
\def\CK{{\mathcal{K}}}

\def\CU{{\mathcal{U}}}

\def\BC{{\mathbf{C}}}

\def\BC{{\mathbb{C}}}
\def\BR{{\mathbb{R}}}
\def\BZ{{\mathbb{Z}}}

\def\wt{\operatorname{wt}\nolimits}

\def\Id{\operatorname{Id}\nolimits}

\def\GL{\operatorname{GL}\nolimits}

\def\Red{\operatorname{Red}\nolimits}
\def\Aut{\operatorname{Aut}\nolimits}

\def\Cat{\operatorname{Cat}\nolimits}

\def\reg{\operatorname{reg}\nolimits}
\def\gen{\operatorname{gen}\nolimits}

\def\im{\operatorname{im}\nolimits}
\def\re{\operatorname{re}\nolimits}

\def\Poin{\operatorname{Poin}\nolimits}

\def\Spec{\operatorname{Spec}\nolimits}

\def\codim{\operatorname{codim}\nolimits}

\def\Disc{\operatorname{Disc}\nolimits}
\def\Flag{\operatorname{Flag}\nolimits}
\def\ie{{\em i.e.}}

\title{Topology of complex reflection arrangements}
\author{David Bessis}
\address{DMA, \'Ecole normale sup\'erieure, 
45 rue d'Ulm, 75230 Paris cedex 05, France}
\email{david dot bessis at ens dot fr}

\begin{document}

\begin{abstract}
Let $V$ be a finite dimensional complex vector space and
$W\subset \GL(V)$ be a finite complex reflection group.
Let $V^{\reg}$ be the complement in $V$ of the reflecting
hyperplanes. A classical conjecture predicts that $V^{\reg}$
is a $K(\pi,1)$ space. When $W$ is a complexified real
reflection group, the conjecture follows from a theorem of Deligne,
\cite{deligne}. Our main result validates the conjecture
for \emph{duality} (or, equivalently, \emph{well-generated})
complex reflection groups. This includes
the complexified real case (but our proof is new) and new cases
not previously known.
We also address a number of questions about $\pi_1(W\cq V^{\reg})$,
the braid group of $W$.
\end{abstract}

\maketitle

Let $V$ be a finite dimensional complex vector space and
$W\subset \GL(V)$ be an irreducible complex reflection group
(all reflection groups considered here are assumed to be finite).
Let $d_1 \leq \dots \leq d_n$ be the degrees of $W$.
Let $d_1^* \geq \dots \geq d_n^*=0$ be the codegrees of $W$.
We say that $W$ is a \emph{duality} group if
$d_i+d_i^*=d_n$ for all $i$ (by analogy with the real case, we say that
 $d_n$ is the
\emph{Coxeter number} of $W$ and we often denote it by $h$).
We say that $W$ is \emph{well-generated} if it may be generated 
by $n$ reflections.
Orlik-Solomon observed, by inspecting the classification of Shephard-Todd,
that
$$W \text{ is a duality group } \Leftrightarrow W \text{ is well-generated}.$$

Let $V^{\reg}$ be the complement in $V$ of the reflecting
hyperplanes.
In the case when $W$ is a type $A$ reflection group,
Fadell and Neuwirth proved in the early 1960's that
$V^{\reg}$ is a $K(\pi,1)$ (this is an elementary use of fibration
exact sequences, \cite{fadell}).
Brieskorn conjectured in 1971, \cite{brieskorn2},
that the $K(\pi,1)$ property holds
when $W$ is a complexified real reflection group.
It is not clear who first stated the conjecture for arbitrary
complex reflection groups. The conjecture may be found in Orlik-Terao's
book:

\begin{conj}[\cite{orlikterao}, p. 163 \& p. 259]
The universal cover of $V^{\reg}$ is contractible.
\end{conj}

The complexified real case (\ie, Brieskorn's conjecture)
 was quickly settled by Deligne, \cite{deligne}.
The rank $2$ case is trivial.
The case of the infinite family was solved in 1983 by Nakamura, \cite{nakamura}
(here again, the monomiality of the group allows an efficient use of
fibrations).
A few other cases immediately follow from the observation by Orlik-Solomon,
\cite{orliksolomon}, that certain orbifolds of non-real
 complex reflection groups
are isomorphic to orbifolds of complexified real reflection groups.

Combining all previously known results, the conjecture remained open for
six exceptional types: $G_{24}$, $G_{27}$, $G_{29}$, $G_{31}$,
$G_{33}$ and $G_{34}$.
Our main result is:
\begin{theo}
\label{theokapiun}
Assume that $W$ is a well-generated complex reflection group. The universal
cover of $V^{\reg}$ is contractible.
\end{theo}

Five of the six open cases are duality groups:
$G_{24}$, $G_{27}$, $G_{29}$, $G_{33}$ and $G_{34}$.
The theorem also applies to the complexified 
real case, for which we obtain a new proof,
not relying on \cite{deligne}.
The only case left is $G_{31}$, an irreducible complex reflection group
of rank $4$ which cannot be generated by less than $5$ reflections.

Among by-products of our construction, we obtain new cases of
several standard conjectures about the braid group of $W$, defined
by  $$B(W) := \pi_1(W \cq V^{\reg}).$$

\begin{theo}
\label{theointrogarside}
Braid groups of well-generated complex reflection groups
are Garside groups.
\end{theo}

In particular, they are torsion-free, have nice solutions to word
and conjugacy problems, are biautomatic, admit finite $K(\pi,1)$
(our construction provides an explicit one), and much more -- see
\cite{dehgar} for a quite complete reference. None of this was known
for the five exceptional groups mentioned above.

\begin{theo}
\label{theocenter}
The center of the braid group of an irreducible well-generated
complex reflection group is cyclic.
\end{theo}

This, as conjectured by Brou\'e-Malle-Rouquier, \cite{bmr}, should
also hold for non-well-generated groups. Again, the cases of
$G_{24}$, $G_{27}$, $G_{29}$, $G_{33}$ and $G_{34}$ are new.
The only case left is $G_{31}$.

For $B(G_{29})$, $B(G_{31})$, $B(G_{33})$ and $B(G_{34})$, no presentations
were known until now, although some conjectures made in \cite{BM} were
supported by strong evidences.

\begin{theo}
\label{theointropres}
The conjectural presentations
for $B(G_{29})$, $B(G_{33})$ and $B(G_{34})$ given in \cite{BM}
are correct.
\end{theo}

\subsection*{General strategy}

The general architecture of our proof is borrowed from
Deligne's original approach but the details are quite different.
Every construction here is an analog
of a construction from \cite{deligne} but relies on different combinatorial
and geometric objects.
Like in \cite{deligne}, one studies a certain braid monoid $M$, whose
structure expresses properties of reduced decompositions in $W$,
and one proves
that it is a lattice for the divisibility order
(this amounts to saying that the monoid is \emph{Garside}). Like in
\cite{deligne}, one uses semi-algebraic geometry to construct
an open covering of the universal cover of $V^{\reg}$, with the property
that non-empty intersections are contractible. This implies that
the universal cover is homotopy equivalent to the nerve of the covering.
Like in \cite{deligne}, one interprets this nerve as a certain flag complex
obtained from $M$. Like in \cite{deligne}, the contractibility of the nerve
follows from the lattice property for $M$.
However, our proof does not use the classical braid monoid, but
a dual braid monoid (\cite{dualmonoid}, \cite{BC}), whose construction
is generalised to all well-generated complex reflection groups.
The construction of the open covering is the most problematic step:
by contrast with the real case, one cannot rely on the notions of \emph{walls}
and \emph{chambers}. The idea here is to work in $W\cq V^{\reg}$ and to use
a generalisation of the Lyashko-Looijenga morphism. This morphism
allows a description of $W\cq V^{\reg}$ by means of a
 ramified covering of a type
$A$ reflection orbifold.
Classical objects like \emph{walls}, \emph{chambers} and \emph{galleries}
can somehow be ``pulled-back'',
via the Lyashko-Looijenga morphism, to give semi-algebraic objects
related to the dual braid monoid.

\subsection*{Lyashko-Looijenga coverings} The
quotient map $\pi:V^{\reg} \twoheadrightarrow W \cq V^{\reg}$ is a regular
covering.
Once a system of basic invariants $(f_1,\dots,f_n)$ is chosen, the quotient
space
$W \cq V^{\reg}$ identifies with the complement in $\BC^n$ 
of an algebraic hypersurface
$\CH$, the \emph{discriminant}, of equation $\Delta\in\BC[X_1,\dots,X_n]$.
If $W$ is an irreducible duality complex reflection group,
it is possible to choose $(f_1,\dots,f_n)$ such that 
$$\Delta= X_n^n + \alpha_2 X_n^{n-2} + \dots 
+ \alpha_n,$$
where $\alpha_2,\dots,\alpha_n\in \BC[X_1,\dots,X_{n-1}]$.
Let $Y := \Spec \BC[X_1,\dots,X_{n-1}]$, together with the natural
map $p:W \cq V^{\reg} \rightarrow Y$.
We have an identification  $W\cq V \simeq \BC^n \stackrel{\sim}{\rightarrow}
\BC\times Y$ sending the orbit $\overline{v}$ of $v\in V$ to
$(f_n(v),p(\overline{v}))$.
The fiber of $p$ over $y\in Y$ is a line $L_y$ which
intersects $\CH$ at $n$ points (counted with multiplicities).
Generically,
the $n$ points are distinct. Let $\CK$ be the \emph{bifurcation locus},
\ie, the algebraic hypersurface of $Y$ consisting of points $y$ such
that the intersection has cardinality $<n$.
Classical results from invariant theory of complex reflection groups
make possible (and very easy) to generalise a
 construction by Looijenga and Lyashko: the map $LL$ (for
``Lyashko-Looijenga'') sending
$y\in Y-\CK$ to the subset $\{x_1,\dots,x_n\} \subseteq
\BC$ such that $p^{-1}(y)\cap \CH  = \{(x_1,y),\dots,(x_n,y)\}$ is a regular
covering of degree $n!h^n/|W|$ of the (centered)
configuration space of $n$ points in $\BC$. In particular,
$Y-\CK$ is a $K(\pi,1)$.
This observation, which is apparently new in the non-real case, already
allows a refinement of our earlier results (\cite{zariski}, \cite{BM})
on presentations for the braid group of $W$.

\subsection*{The dual braid monoid}
When $W$ is complexified real, a \emph{dual braid monoid} was constructed in
\cite{dualmonoid} (generalising the construction of Birman-Ko-Lee, \cite{bkl};
similar partial results were
independently obtained by Brady and Brady-Watt, \cite{bw}).
The construction was later generalised in \cite{BC}
to the complex reflection group $G(e,e,n)$.
Assume $W$ is a
well-generated complex reflection group generated by reflections of
order $2$. Let $R$ be the set of all reflections in $W$.
The idea is that the pair $(W,R)$ has some ``Coxeter-like'' features.
Instead of looking at relations of the type
 $$\underbrace{sts...}_{m_{s,t}}=\underbrace{tst...}_{m_{s,t}},$$ one
considers relations of the type $$st=tu$$ where $s,t,u\in R$.
Let $\mathcal{S}$ be the set of all relations of this type holding in $W$.
In general, $B(W) \not \simeq \left< R \left|  \mathcal{S} \right.
\right>$, but it is possible to find natural subsets $R_c \subseteq R$ and
$\mathcal{S}_c \subseteq
\mathcal{S}$ such that $B(W) 
\simeq \left< R_c \left|  \mathcal{S}_c \right.
\right>$
(if $W$ is complexified real, $R_c=R$).
The elements of $\mathcal{S}_c$ are called \emph{dual braid relations}.
The choices of $R_c$ and $\mathcal{S}_c$ are  
natural once a \emph{Coxeter element} $c$ has been chosen
(the notion of \emph{Coxeter element} generalises to the non-real
groups in terms of Springer's theory of regular elements, \cite{springer}).
Since the relations are positive, one may view the presentation as
a monoid presentation, defining a monoid $M(W)$.
The crucial property of this monoid is that it is a lattice for the
divisibility order or, in other words, a \emph{Garside monoid}.
Following Deligne, Bestvina, T. Brady and Charney-Meier-Whittlesey,
the Garside structure provides a convenient simplicial Eilenberg-McLane
$K(B(W),1)$ space (\cite{deligne}, \cite{bestvina}, \cite{brady}, \cite{cmw}).
The earlier results on the dual braid monoid are improved here in two
directions:
\begin{itemize}
\item The construction is generalised to the few exceptional cases 
($G_{24}$, $G_{27}$, $G_{29}$,
$G_{33}$ and $G_{34}$) not
covered by \cite{dualmonoid} and \cite{BC}.
\item A new geometric interpretation is given, via the Lyashko-Looijenga
covering. This interpretation is different for the one given in
\cite[Section 4]{dualmonoid}.
\end{itemize}

The second point is the most important. It relies (so far) 
on a counting argument,
following and extending a property which, for the complexified real case,
was conjectured by Looijenga and proved
in a letter from Deligne to Looijenga, \cite{deligneletter}.

\subsection*{Tunnels}
The classical theory of real reflection groups combines a
``combinatorial'' theory (Coxeter systems) and a ``geometric''
theory (expressed in the language, invented by Tits, of
walls, chambers, galleries, buildings...)
We expect the dual braid monoid approach to eventually provide
effective substitutes for much of this classical theory.
We use here a notion of \emph{tunnel}, which is a rudimentary
geometric object replacing
the classical notion of \emph{minimal gallery between two chambers}.
An important difference with the classical geometric language is
that tunnels are naturally visualised in $W\cq V$ (instead
of $V$).
A tunnel $T$ 
is a path in $W\cq V^{\reg}$ drawn inside a single
line $L_y$ (for some $y\in Y$) and with
constant imaginary part. It represents an element $b_T$ of the dual braid
monoid $M$. An element of $M$ is \emph{simple} if it is represented by
a tunnel. This notion coincides with the notion of simple element 
associated with the Garside structure.
In the classical approach, for any chamber $\mathcal{C}$, there are as
many equivalence classes of 
minimal galleries starting at $\mathcal{C}$ as simple elements (this number
is $|W|$). Here the situation is different: in a given $L_y$, not all
simples are represented. The simples represented in different $L_y$'s
may be compared thanks to a huge ``basepoint'' $\CU$ which
is both dense in $W\cq V^{\reg}$ and contractible.

\subsection*{Non-crossing partitions}
In the classical cases $A_n$, $B_n$, $D_n$ and, more generally,
$G(e,e,n)$, the structure of $M(W)$ is understood in terms
of suitable notions of \emph{non-crossing partitions}
(\cite{athanasiadisreiner}, \cite{BC}, \cite{BDM}). The dual braid
monoid of an irreducible well-generated complex reflection $W$ gives rise
to a lattice of \emph{generalised non-crossing partitions},
whose cardinal is the \emph{generalised Catalan number}
$$\Cat(W) := \prod_{i=1}^n \frac{d_i+d_n}{d_i}.$$
(the term ``partition'' should not be taken too seriously: the elements
of the lattice no longer have a natural interpretation as partitions.)
It is likely that this combinatorial object has some representation-theoretic
interpretation. If $W$ is not well-generated, $\Cat(W)$ may fail to be
an integer.

\subsection*{Acknowledgments} Work on this project started during
stimulating visits to KIAS (Seoul) and RIMS (Kyoto), in the summer of 2003,
when I realised that the dual braid monoid construction
extends to duality
groups. I thank Sang Jin Lee and Kyoji Saito for their hospitality and
for their interests in discussing these subjects.
After this initial progress, I remained stucked for many months, trying
to construct the open covering using affine geometry in $V^{\reg}$
(refining \cite[Section 4]{dualmonoid}).
Two observations helped me decide to work in the quotient. 
First, Kyoji Saito pointed out that the starting point for the
construction of the \emph{flat structure} (or \emph{Frobenius manifold
structure}, see \cite{saitoorbifold}) on real reflection orbifolds 
was precisely 
the duality between degrees and codegrees.
The intuition that the flat structure has something to do with the 
$K(\pi,1)$ property is explicitly mentioned as a motivation for
 \cite{saitoorbifold} (see also \cite[p. 7]{saitopolyhedra}).
The second useful discussion was with Fr\'ed\'eric Chapoton, who
pointed out the numerological coincidence, in the Coxeter case, between the
degree of the Lyashko-Looijenga covering and the number of
maximal chains in the lattice of
non-crossing partitions.

\subsection*{Warning}
This text was initially planned to be a research announcement.
Though it grew much beyond any reasonable size for a research announcement,
it is not yet in final form. 
It still misses some details, and several arguments are not optimal.

\section*{Notations}

We save the letter $i$ for indexing purposes and denote
by $\sqrt{-1}$ a complex square root of $-1$ fixed once for all.
If $n$ is a positive integer, we denote by $\zeta_n$ the standard $n$-th
root of unity $\exp(2\sqrt{-1}\pi /n)$.

Many objects depend on a complex reflection group $W$, e.g., the braid
group $B(W)$. We often drop the explicit mention of $W$, and write
$B$ for $B(W)$. When $n$ is an integer, we denote by $B_n$ the braid
group on $n$ strings, together with its standard generating set
$\bm{\sigma}_1,\dots,\bm{\sigma}_{n-1}$; it is isomorphic to the braid group
of $\mathfrak{S}_n$ in its permutation reflection representation
(see Section \ref{sectionconfig}). The groups $B$ and $B_n$ appear
simultaneously and should not be confused.
 
\section*{Topological conventions}

Let $E$ be a topological space. Let $\gamma$ be a path in $E$, \ie,
a continuous map $[0,1]\rightarrow E$. We say that $\gamma$ is a \emph{path
from $\gamma(0)$ to $\gamma(1)$} -- or that $\gamma(0)$ is the \emph{source}
and $\gamma(1)$ the \emph{target}. The concatenation rule is as follows:
if $\gamma,\gamma'$ are paths such
that $\gamma(1)=\gamma'(0)$,
we define by $\gamma\circ  \gamma'$ the path mapping
$t\leq 1/2$ to $\gamma(2t)$ and $t\geq 1/2$ to $\gamma'(2t-1)$.

A path $\gamma$ from $e$ to $e'$ induces an isomorphism $\phi_{\gamma}:
\pi_1(E,e) \stackrel{\sim}{\rightarrow} \pi_1(E,e')$.

Let $X$ be a subspace of $E$.
Assume that $X$ is simply-connected (in all applications, $X$ will
actually be contractible).
Then for all $x,x'\in X$
there is a unique homotopy class of paths from $x$ to $x'$ drawn in $X$,
thus a natural isomorphism
$\phi_{X,x,x'}: \pi_1(E,x) \stackrel{\sim}{\rightarrow} \pi_1(E,x')$.
By unicity, $\phi_X:=(\phi_{X,x,x'})_{x,x'\in X}$ is a transitive system
of isomorphisms between the $(\pi_1(E,x))_{x\in X}$. We denote
by $$\pi_1(E,X)$$ the transitive limit.
Practically speaking, $\pi_1(E,X)$ should be thought of
as any $\pi_1(E,x)$, for some $x\in X$, together with an unambiguous
receipe thanks to which any path in $E$
from any $x'\in X$ to any $x''\in X$ represents a unique element
of $\pi_1(E,x)$ -- moreover, one may
forget about which $x\in X$ was chosen and change it at our convenience.

Assume that $E$ is locally simply-connected.
To get a model $\widetilde{E}_x$ 
of the universal cover $\widetilde{E}$, one chooses a
basepoint $x\in E$. 
A path $\gamma$ from $x$ to $e\in E$ represents a point
$\overline{\gamma}\in
\widetilde{E}_x$
in the fiber over $e$. Two paths represent the same point if and only if
they are homotopic (with fixed endpoints).
The topology of $\widetilde{E}_x$ around $\overline{\gamma}$
is described by concatenation at the target
with paths starting at $e$ and drawn
in simply-connected neighbourhoods of $e$.
The fundamental group $\pi_1(E,x)$ acts on the left on $\widetilde{E}_x$ by
concatenation at the source. More generally, a path from $x$ to $x'$
induces an isomorphism $\widetilde{E}_x\stackrel{\sim}{\rightarrow}
\widetilde{E}_{x'}$. Here again, if $X$ is simply-connected,
we have a transitive system of isomorphisms between
$(\widetilde{E}_x)_{x\in X}$ whose limit is a model $\widetilde{E}_X$ of
the universal cover. Practically speaking, any path from any
$x\in X$ to any $e\in E$ unambiguously represents a 
point in $\widetilde{E}_X$ (in the fiber over $e$).
Moreover, $\pi_1(E,X)$ acts on the left on $\widetilde{E}_X$.

\section{Complex reflection groups, discriminants, braid groups}

Let $V$ be a $\BC$-vector space of finite dimension $n$.
A \emph{complex reflection group} in $\GL(V)$ is a (finite)
subgroup generated by \emph{complex reflection}, elements whose fixed
subspace is a hyperplane.
Let $W\subset \GL(V)$ be a complex reflection group.
A \emph{system of basic invariants} is an $n$-tuple  $f=(f_1,\dots,f_n)$
of homogeneous $W$-invariant polynomial functions on $V$ generating
the algebra $\BC[V]^W$ of $W$-invariant polynomial functions.
A classical theorem of Shephard-Todd, \cite{shephard}, asserts
that such tuples exist, and that they consist of algebraically independent
terms.
Set $d_i:=\deg f_i$; these numbers are the \emph{degrees} of $W$.
Up to reordering, we may assume that
$d_1\leq d_2 \leq \dots \leq d_n$. The sequence $(d_1,\dots,d_n)$ is
then independent of the choice of $f$.

Choosing a system of basic invariants $f$ amounts to choosing
a graded algebra isomorphism
$\BC[V]^W \simeq\BC[X_1,\dots,X_n], f_i \mapsto X_i$,
where the indeterminate
$X_i$ is declared homogeneous with degree $d_i$. Geometrically,
this isomorphism identifies the categorical quotient $W\cq V$
with the affine space $\BC^n$.

Further features of the invariant theory of complex reflection
groups involve invariant vector fields and invariant differential forms
on $V$. 

\begin{theo}[\cite{orlikterao}, Lemma 6.48]
The $\BC[V]^W$-modules $(\BC[V]\otimes V)^W$ and $(\BC[V]\otimes V^*)^W$ are
free of rank $n$.
\end{theo}

If $f=(f_1,\dots,f_n)$ is a system of basic invariants, $df:=(df_1,\dots,df_n)$
is a $\BC[V]^W$-basis for $(\BC[V]\otimes V^*)^W$.
Being homogeneous, the module $(\BC[V]\otimes V)^W$ admits an homogeneous
basis. 

\begin{defi}
A system of basic derivations is an homogeneous $\BC[V]^W$-basis
$\xi=(\xi_1,\dots,\xi_n)$ of $(\BC[V]\otimes V)^W$,
with $\deg(\xi_1) \geq \deg(\xi_2) \geq \dots \geq \deg(\xi_n)$.

The sequence $(d_1^*,\dots,d_n^*):=
(\deg(\xi_1),\dots,\deg(\xi_n))$ (which does not depend on the choice
of $\xi$) is the sequence of \emph{codegrees} of $W$.
\end{defi}

Note that, as in \cite{zariski}, we label codegrees in decreasing order,
which is slightly unusual.
When $W$ is a complexified real reflection group, we have $V\simeq V^*$
as $W$-modules, thus $d_i^* = d_{n-i+1}-2$ for all $i$. This relation
is specific to the real situation and is not relevant here.

The Euler vector field on $V$ is invariant and of degree $0$.
Thus $d_n^*=0$.

Invariant vector fields define vector fields on the quotient variety.
Let $f$ be a system of basic invariants and $\xi$ be a system of
basic derivations.
For $j\in\{1,\dots,n\}$, the vector field $\xi_j$ defines a vector
fiels $\overline{\xi_j}$ on $W \cq V$. Since $\frac{\partial}{\partial f_1},
\dots,\frac{\partial}{\partial f_n}$ is a $\BC[V]^W$-basis of the module
of polynomial vector fields on $W \cq V$, we have
$$\overline{\xi_j} = \sum_{i=1}^n m_{i,j} \frac{\partial}{\partial f_i},$$
where the $m_{i,j}$ are uniquely defined elements of $\BC[V]^W$.

\begin{defi}
The discriminant 
matrix of $W$ (with respect to $f$ and $\xi$) is $M:=(m_{i,j})_{i,j}$.
\end{defi}

To recover $\xi_j$ from $\overline{\xi_j}$, replace
$\frac{\partial}{\partial f_i}$ by $\sum_{k=1}^n \frac{\partial f_i}{\partial
Y_k}\frac{\partial}{\partial Y_i}$.
By weighted homogeneity, one has:
\begin{lemma}
\label{wtmij}
For all $i,j$, $\wt(m_{i,j})= d_i + d_j^*$.
\end{lemma}

The vector space $V$ decomposes as a direct
sum $\bigoplus_i V_i$ of irreducible representations of $W$.
Denote by $W_i$ the irreducible reflection group in $\GL(V_i)$
generated by (the restriction of) the reflections in $W$
whose hyperplanes contain $\bigoplus_{j\neq i}V_j$.
We have $W \simeq \prod_{i} W_i$.
Viewing $W$ and the $W_i$'s as
\emph{reflection groups}, \ie, groups endowed with a reflection representation,
it is natural to actually write $W= \bigoplus_iW_i$.

We denote by $\CA$ the arrangement of $W$, \ie, the
set of reflecting hyperplanes of reflections in $W$;
it is in bijection with the disjoint union $\sqcup_i \CA_i$, where
$\CA_i$ is the arrangement of $W_i$.
We set $$V^{\reg}:=V-\bigcup_{H\in \CA}H.$$
Denote by $p$ the quotient map $V\twoheadrightarrow W \cq V$.
Choose a basepoint $v_0\in V^{\reg}$.

\begin{defi}[\cite{bmr}] The \emph{braid group} of $W$
is $B(W):=\pi_1(W\cq V^{\reg},p(v_0))$.
\end{defi}

To write explicit equations, one chooses a system of 
basic invariants $f$. The \emph{discriminant}
$\Delta(W,f)\in\BC[X_1,\dots,X_n]$ is
the reduced equation of $p(\bigcup_{H\in \CA}H)$, via the 
identification $\BC[V]^W \simeq\BC[X_1,\dots,X_n]$.

One easily sees that $B(W)\simeq \prod_{i} B(W_i)$. More generally,
all objects studied here behave ``semi-simply'', and we may restrict
our attention to irreducible complex reflection groups.

Since $B(W)$ is the fundamental group of the complement of an 
algebraic hypersurface, it is generated by particular elements
called \emph{generators-of-the-monodromy} or \emph{meridiens}
(see, for example, \cite{bmr} or \cite{zariski}). They map
to reflections under the natural epimorphism $B(W)\rightarrow W$.
The diagrams given in \cite{bmr} symbolise presentations whose generators
are generators-of-the-monodromy (except for the six exceptional types
for which no presentation was known).
\begin{defi}
The generators-of-the-monodromy of $B(W)$ are called \emph{braid reflections}.
\end{defi}
This terminology was suggested by Brou\'e. It is actually tempting
to simply call them \emph{reflections}: since they generate $B(W)$, the
braid group appears to be some sort of ``reflection group''.

Another natural feature of $B(W)$ is the existence of a natural
length function, which is the unique group morphim
$$l:B(W) \rightarrow \BZ$$
such that, for all braid reflection $s\in B(W)$, $l(s)=1$.

Consider the 
intersection lattice $\CL(\CA):=\{ \bigcap_{H\in A} H | A\subseteq \CA\}$.
Elements of $\CL(\CA)$ are called \emph{flats}.
It is standard to endow $\CL(\CA)$ with the reversed-inclusion partial
ordering:
$$ \forall L,L' \in \CL(\CA), L \leq L' :\Leftrightarrow L \supseteq L'.$$

For $L \in \CL(\CA)$, we denote by $L^0$ the complement in $L$ of
the flats strictly included in $L$.
The $(L^0)_{L \in \CL(\CA)}$ form a stratification $\CS$ of $V$.
We consider the partial ordering on $\CS$ defined by 
$L^0 \leq L'^0 :\Leftrightarrow L \leq L'$. This is a degeneracy relation:
$$\forall L\in\CL(\CA), \overline{L^0} = L= \bigcup_{S \in \CS, L^0 \leq S} S.$$

Since $W$ acts on $\CA$, it acts on $\CL(\CA)$ and we obtain 
a quotient stratification $\overline{\CS}$ of $W\cq V$ called
\emph{discriminant stratification}.

\begin{prop}[\cite{orlikterao}, Corollary 6.114]
\label{span}
Let $v\in V$. The vectors $\xi_1(v),\dots,\xi_n(v)$ span the
tangent space to the stratum of $\CS$ containing $v$.
The vectors $\overline{\xi}_1(v),\dots,\overline{\xi}_n(v)$ span the
tangent space to the stratum of $\overline{\CS}$ containing $\overline{v}$.
\end{prop}

\section{Well-generated complex $2$-reflection groups}

Irreducible complex reflection groups
were classified fifty years ago by Shephard and Todd, \cite{shephard}.
There is an infinite family $G(de,e,n)$, where $d,e,n$ are positive
integers, and $34$ exceptions $G_4,\dots,G_{37}$.
Let us distinguish three subclasses of complex reflection groups:
\begin{itemize}
\item \emph{real reflection groups}, obtained by scalar extension from
reflection groups of real vector spaces;
\item \emph{$2$-reflection groups},
 generated by reflections of order $2$;
\item \emph{well-generated reflection groups}, \ie,
irreducible complex reflection groups generated by $n$ reflections,
where $n$ is the dimension, and more generally direct sums of such irreducible
groups.
\end{itemize}
Real reflection groups are both $2$-reflection groups and well-generated.
For non-real groups, any combination of the other two properties may hold.

As far as the $K(\pi,1)$ conjecture and properties of 
braid groups are concerned, it is enough to restrict
one's attention to $2$-reflection groups, thanks to the following theorem.

\begin{theo}
\label{shephard}
For any complex reflection group $W\subset \GL(V)$, one may find
$V'$ a complex vector space, $W'\subset \GL(V')$ a complex
$2$-reflection group and $f$ (resp. $f'$) a system
of basic invariants 
for $W$ (resp. $W'$) such that $\Delta(W,f)=\Delta(W',f')$.
In particular, $W\cq V^{\reg}\simeq W'\cq V'^{\reg}$ and
$B(W)\simeq B(W')$.
\end{theo}

\begin{proof}
This theorem may be observed on the classification and
was certainly known to experts but remained unexplained.
The recent work of Couwenberg-Heckman-Looijenga, \cite{chl}, may be adapted
to provide a direct
argument, which is sketched below.

All references and notations are from \emph{loc. cit.}.
Assume that $W$ is not a $2$-reflection group.
For each $H\in \CA$, let $e_H$ be the order of the pointwise stabiliser
$W_H$ and set $\kappa_H:=1 - e_H/2$. Consider the Dunkl connection
$\nabla$ with connection form $\sum_{H\in \CA} \omega_H\otimes \kappa_H
\pi_H$, as in Example 2.5. Since $e_H\geq 2$, we have $\kappa_H\leq 0$.
In particular, $\kappa_0=1/n\sum_{H\in\CA}\kappa_H \leq 0$ (Lemma 2.13)
and we are in the situation of \emph{loc. cit.} Section 5. 
In many cases, this suffices to conclude.
The problem is that, even though at least
some $e_H$ have to be $>2$, it is possible
that $\CA$ contains several orbits, some of them with $e_H=2$. To handle
this, one has to enlarge the ``Schwarz symmetry group'' of
\emph{loc. cit.}, Section 4.
\end{proof}

The importance of the distinction between well-generated and
``badly-generated'' groups was first pointed out by Orlik-Solomon,
who observed in \cite{orso} 
a coincidence with invariant-theoretical aspects.
\begin{theo}
\label{monic}
Let $W$ be an irreducible complex reflection group. The following
assertions are equivalent:
\begin{itemize}
\item[(i)] $W$ is well-generated.
\item[(ii)] For all $i\in \{1,\dots,n\}$,
$d_i+d_i^* = d_n$.
\item[(iii)] For all $i\in \{1,\dots,n\}$, $d_i+d_i^* \leq d_n$.
\item[(iv)] For any system of basic invariants $f$, there exists
a system
of basic derivations $\xi$ such that the discriminant matrix decomposes as
$M= M_0 + X_n M_1$, where $M_0,M_1$ are matrices with coefficients
in $\BC[X_1,\dots,X_{n-1}]$ and $M_1$ is lower triangular with non-zero
scalars
on the diagonal.
\item[(v)] For any system of basic invariants $f$, we have
$\frac{\partial^n\Delta(W,f)}{(\partial X_n)^n}\in \BC^{\times}$
(in other words,
$\Delta(W,f)$, viewed as a polynomial in $X_n$ with coefficients
in $\BC[X_1,\dots,X_{n-1}]$, is monic of degree $n$).
\end{itemize}
\end{theo}

The matrix $M_1$ from assertion (iv) is an analog of the matrix $J^*$
from \cite{saitoorbifold}, p. 10. Assertion (iv) itself generalises the
non-degeneracy argument for $J^*$, which is an 
important piece of the construction of Saito's ``flat structure''.

\begin{proof}
(i) $\Rightarrow$ (ii) was observed in \cite{orso}
inspecting the classification. We still have no good explanation.

(ii) $\Rightarrow$ (iii) is trivial.

(iii) $\Rightarrow$ (v).
Let $h:=d_n$. A first step is to observe that,
under assumption (iii), $h$ is a regular number, in the sense of Springer,
\cite{springer}. Indeed, since $d_n^*$ is the only codegree equal to $0$,
the condition $d_i+d_i^*\leq h$ implies, for $i=1,\dots,n-1$,
 both $0<d_i< h$ and $0<d_i^*<h$; in particular, except when $i=n$, $h$ does
divides neither $d_i$ nor $d_i^*$.
A criterion first noticed by Lehrer-Springer
and proved abstractly by Lehrer-Michel (\cite{lehrermichel}) says that if
a number $d$ divides exactly as many degrees as codegrees, then it is regular.
Since $h$ is regular and divides only one degree, we may use
\cite{zariski}, Lemma 1.6 (ii) to obtain assertion (v):
the discriminant is $X_n$-monic, and by weighted-homogeneity it must be
of degree $n$.

(iii) $\Rightarrow$ (iv) is a refinement of the previous discussion.
Each entry $m_{i,j}$ of the matrix $M$ is weighted-homogeneous of
total weight $d_i+d_j^* \leq d_n+d_1^*= 2h - d_1 < 2h$; since
$X_n$ has weight $h$, $\deg_{X_n}m_{i,j} \leq 1$. This explains
the decomposition $M = M_0 + X_n M_1$, where $M_0$ and $M_1$ have
coefficients in $\BC[Y]$.
If $i<j$ and $d_i< d_j$, then $d_i+d_j^* < d_j + d_j^* = h$, thus $\deg_{X_n}
m_{i,j} = 0$. The matrix $M_0$ is \emph{weakly lower triangular},
\ie, lower triangular except that
there could be non-zero terms above the diagonal in square diagonal
blocks corresponding to successive equal degrees (successive degrees
may indeed be equal, as in the example of type $D_4$, where the degrees
are $2,4,4,6$).

Let $i_0< j_0$ such that $d_{i_0} =d_{i_0+1} = \dots =d_{j_0}$
(looking at the classification, one may observe
 that this forces $j=i+1$; this observation is not used in the argument
below).
For all $i,j \in
\{i_0,\dots,j_0\}$, we have $d_i+d_j^* = h$. By weighted homogeneity,
this implies that the corresponding square block of $M_1$ consists of scalars.
The basic derivations $\xi_{i_0},\dots,\xi_{j_0}$ all have the same degrees,
thus one is allowed to perform Gaussian elimination on the corresponding
columns of $M$. Thus, up to replacing $\xi$ by another system of basic
derivations $\xi'$, we may assume that $M_1$ is lower triangular.

The diagonal terms of $M_1$ must be scalars, once again by weighted homogeneity.
Assuming (iii), we already know that (v) holds. The determinant
of $M$ is $\Delta(X,f)$; (v) implies that the coefficient of $X_n^n$
is non-zero.
This coefficient is the product of the diagonal terms of $M_1$. We have
proved (iii).

(iv) $\Rightarrow$ (v) is trivial.

(v) $\Rightarrow$ (i) follows from the main result in \cite{zariski}.
\end{proof}

The following notion was considered in \cite{saito} for real
reflection groups. 

\begin{defi}
A system of basic derivations is \emph{flat} (with respect to 
$f$) if
the discriminant matrix may be written $M=M_0+X_n\Id$,
where $M_0$ is a matrix with coefficients in $\BC[X_1,\dots,X_{n-1}]$.
\end{defi}

Theorem \ref{monic} has the following corollary:

\begin{coro}
\label{flat}
Let $W$ be an well-generated irreducible reflection group.
The flatness of a system of basic derivations does not depend
from the choice of $f$. There exists flat systems of basic derivations.
\end{coro}

\begin{proof}
The independence follows from the fact that $d_{n-1}< d_n$,
which itself follows from characterisation (ii) of
well-generated groups.

Existence:
choose $f$ and $\xi$,
and write $M=M_0+X_nM_1$, as in characterisation $(iv)$.
The matrix $M_1$ is invertible
in $\GL_n(\BC[X_1,\dots,X_{n-1}])$.
The matrix $M_1^{-1} M = M_1^{-1}M_0 + X_n \Id$ represents
a flat system of basic derivations (weighted homogeneity is preserved
by the Gaussian elimination procedure).
\end{proof}

Contrary to what happens with real reflection groups, we may
not use the identification $V\simeq V^*$ to obtain ``flat system
of basic invariants''.

Irreducible groups which are not well-generated may always
be generated by $\dim V + 1$ reflections. This fact has been observed
long ago, by case-by-case inspection, but no
general argument is known. (In some sense, these badly-generated groups
should be thought of as affine groups).
The simplest example of a non-well-generated group is the group $G(4,2,2)$,
generated by
$$\begin{pmatrix} -1 & 0 \\ 0 & 1 \end{pmatrix}, \;
\begin{pmatrix} 0 & 1 \\ 1 & 0 \end{pmatrix}, \;
\begin{pmatrix} 0 & -i \\ i & 0 \end{pmatrix}.$$

In the sequel, several arguments are case-by-case. As explained above,
one has to consider only well-generated irreducible complex
$2$-reflection groups. Their list consists of:
\begin{itemize}
\item the Coxeter types $A_n$, $B_n$, $D_n$, $E_6$, $E_7$, $E_8$, $F_4$,
$H_3$, $H_4$ and $I_2(e)$,
\item the infinite family $G(e,e,n)$, where $e$ and $n$ are integers
$\geq 2$, which interpolates the Coxeter series $D_n$ (when $e=2$) and
$I_2(e)$ (when $n=2$),
\item five non-real exceptional types $G_{24}$, $G_{27}$, $G_{29}$, $G_{33}$
and $G_{34}$.
\end{itemize}

\section{Symmetric groups, configurations spaces and classical braid groups}
\label{sectionconfig}

This section is only included to clarify some terminology and notations.
Everything here is classical and elementary.

Let $n$ be a positive integer.
The symmetric group $\mathfrak{S}_n$ may be viewed as reflection
group, acting on $\BC^n$ by permuting the canonical basis.
This representation is not irreducible. Let $H$ be the hyperplane
of equation $\sum_{i=1}^n X_i$ (where $X_1,\dots,X_n$ is the dual canonical
basis of $\BC^n$). It is preserved by $\mathfrak{S}_n$, which acts
on it as an irreducible complex reflection group.

We have
$\BC[X_1,\dots,X_n]^{\mathfrak{S}_n} = \BC[\sigma_1,\dots,\sigma_n]$ and
$\BC[H]^{\mathfrak{S}_n} = \BC[\sigma_1,\dots,\sigma_n]/\sigma_1$,
where $\sigma_1,\dots,\sigma_n$ are the elementary symmetric functions
on $X_1,\dots,X_n$. Set
$E_n':= \mathfrak{S}_n \cq \BC^n = \Spec \BC[\sigma_1,\dots,\sigma_n]$ and 
$E_n := \mathfrak{S}_n \cq H = \Spec \BC[\sigma_1,\dots,\sigma_n]/\sigma_1$.
These spaces have more convenient descriptions in terms of multisets.

Recall that a \emph{multiset}
is a set $S$ (the \emph{support} of the multiset) together with
a map $m:S \rightarrow \BZ_{\geq 1}$ (the \emph{multiplicity}).
The \emph{cardinal} of such a multiset is $\sum_{s\in S} m(s)$
(it lies in $\BZ_{\geq 0} \cup \{\infty\}$).
If $(S,m)$ and $(S',m')$ are two multisets and if $S,S'$ are subsets
of a common ambient set, then we may define a \emph{multiset (disjoint)
union}
$(S,m)\cup (S',m')$, whose support is $S\cup S'$ and whose multiplicity
is $m+m'$ (where $m$, resp. $m'$, is extended by $0$ outside $S$, resp
$S'$).

Let $(x_1,\dots,x_n)\in\BC^n$. The associated $\mathfrak{S}_n$-orbit
is uniquely determined by the multiset
$\{\{x_1,\dots,x_n\}\}:=\bigcup_{i=1}^n (\{x_i\},1)$.
This identifies $E'_n$ with the set of multisets of cardinal $n$
with support in $\BC$ (such multisets are called \emph{configurations
of $n$ points in $\BC$}).
The subvariety $E_n$, defined by $\sigma_1=0$, consists of
configurations
 $\{\{x_1,\dots,x_n\}\}$ with $\sum_{i=1}^nx_i = 0$ (such configurations
are \emph{centered}).
The natural inclusion $E_n \subseteq E'_n$ admits the retraction $\rho$
defined by $$\rho(\{\{x_1,\dots,x_n\}\}):=
\{\{x_1 - \sum_{i=1}^nx_i/n,\dots,x_n - \sum_{i=1}^nx_i/n\}\}.$$
Algebraically, this corresponds to the identification of $\BC[\sigma_1,\dots,
\sigma_n]/\sigma_1$ with $\BC[\sigma_2,\dots,\sigma_n]$.

We find it convenient to use
configurations in $E'_n$ to represent
elements of $E_n$, implicitly working through $\rho$. 
E.g., in the proof of Proposition \ref{YCcontractible},
it makes sense to describe
a deformation retraction of a subspace of $E_n$ to a point in terms
of arbitrary configurations because the construction, which only
implies the relative values of the $x_i$'s,
is compatible with $\rho$. We adopt this viewpoint from now on,
without further justifications (compatibility will always be obvious).

Consider the lexicographic total ordering of $\BC$:
if $z,z'\in \BC$, we set
$$z \leq z' :\Leftrightarrow \left\{ \begin{matrix}
\re(z) < \re(z') & \text{or} \\
\re(z)=\re(z') \text{ and } \im(z) \leq \im(z'). \end{matrix} \right.$$

\begin{defi}
The \emph{ordered support} of an element of $E_n$ is the unique sequence
$(x_1,\dots,x_k)$ such that $\{x_1,\dots,x_k\}$ is the support and
$x_1 < x_2 < \dots < x_k$.
\end{defi}

We may uniquely represent an element of
$E_n$ by its ordered support $(x_1,\dots,x_k)$ and the sequence
$(n_1,\dots,n_k)$ of multiplicities at $x_1,\dots,x_k$.

The regular orbit space $E_n^{\reg}:=\mathfrak{S}_n\cq H^{\reg}$
consists of those multisets whose support
has cardinal $n$ (or, equivalently, whose multiplicity 
is constantly equal to $1$).
More generally, the strata of the discriminant stratification of $E_n$
are indexed by partitions of $n$: the stratum $S_{\lambda}$
associated with a
partition $\lambda=(\lambda_1,\lambda_2,\dots,\lambda_k)$, where
the $\lambda_j$'s are integers with
$\lambda_1\geq \lambda_2 \geq \dots \geq \lambda_k > 0$
and  $\sum_{j=1}^k\lambda_j=n$,
consists of configurations whose supports has cardinal $k$
and whose multiplicity functions take the
values $\lambda_1,\dots,\lambda_k$ (with multiplicities).

The braid group $B_n$ associated with $\mathfrak{S}_n$ is the usual
braid group on $n$ strings. We need to be more precise about our
choice of basepoint. For this purpose, we define
$$E_n^{\gen}$$ as
the subset of $E_n^{\reg}$ consisting of configurations of $n$
points with distinct real parts (this is the first in a series of definitions
of semi-algebraic nature). It is clear that:
\begin{lemma}
\label{lemmaEngen}
$E_n^{\gen}$
is contractible.
\end{lemma}

Using our topological conventions, we set $$B_n:=\pi_1(E_n,E_n^{\gen}).$$
This group admits a \emph{standard generating set} (the one
considered by Artin), consisting of
braid reflections $\bm{\sigma}_1,\dots,\bm{\sigma}_{n-1}$ 
defined as follows (we used bold fonts to avoid confusion with the elementary
symmetric functions).
Let $(x_1,\dots,x_n)$ be the ordered support of a point in $E_n^{\gen}$.
Then $\bm{\sigma}_i$
is represented by the following motion of the support:
$$
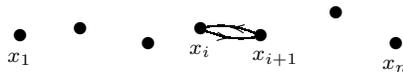
\xy
(-10,0)="1", (-2,1)="2", (7,-1)="3", 
(14,1)="4", (22,0)="5", (32,3)="6", (40,-1)="7",
"1"*{\bullet},"2"*{\bullet},"3"*{\bullet},"4"*{\bullet},
"5"*{\bullet},"6"*{\bullet},"7"*{\bullet},
(14,-2)*{_{x_i}},(24,-3)*{_{x_{i+1}}},
(-10,-3)*{_{x_1}},(40,-4)*{_{x_{n}}},
"4";"5" **\crv{(18,-1)}
 ?(.5)*\dir{>},
"5";"4" **\crv{(18,2)}
 ?(.5)*\dir{>}
\endxy $$

Artin's presentation for $B_n$ is
$$B_n  = \left< \bm{\sigma}_1,\dots,\bm{\sigma}_{n-1}
 \left| \bm{\sigma}_i\bm{\sigma}_{i+1}
\bm{\sigma}_i = \bm{\sigma}_{i+1}\bm{\sigma}_i
\bm{\sigma}_{i+1}, \bm{\sigma}_i\bm{\sigma}_j = 
\bm{\sigma}_j\bm{\sigma}_i \text{ if } |i-j| > 1
 \right. \right>.$$
The following definition requires a compatibility condition which
is a classical elementary consequence of the above presentation.

\begin{defi}
Let $G$ be a group. The \emph{(right)
Hurwitz action of $B_n$ on $G^n$}
is defined by
$$(g_1,\dots,g_{i-1},g_i,g_{i+1},g_{i+2},\dots,g_n) \cdot \bm{\sigma}_i
:= (g_1,\dots,g_{i-1},g_{i+1},g_{i+1}^{-1}g_ig_{i+1},g_{i+2},\dots,g_n),$$
for all $(g_1,\dots,g_n)\in G^n$ and all $i\in \{1,\dots,n-1\}$.
\end{defi}

This action preserves the fibers of the product
map $G^n\rightarrow G,(g_1,\dots,g_n)
\mapsto g_1\dots g_n$.

\section{Lyashko-Looijenga coverings}
\label{sectionlyashkolooijenga}

Let $W$ be an irreducible well-generated complex reflection group, together with
a system of basic invariants $f$ and a flat system of basic
derivations $\xi$ with discriminant matrix $M_0+X_n\Id$.
Expanding the determinant, we observe that
$$\Delta_f = \det ( M_0+X_n \Id)= 
 X_n^n +\alpha_2 X_n^{n-2} + 
\alpha_3 X_n^{n-3} + \dots + \alpha_n,$$
where $\alpha_i\in \BC[Y]$.
Since $\Delta_f$ is weighted homogeneous of total weight $nh$ for the system
of weights $\wt(X_i)=d_i$, each $\alpha_i$ is weighted homegeneous of
weight $ih$. 

\begin{defi}
The (generalised) \emph{Lyashko-Looijenga morphism} is the
morphism $LL$ from $Y = \Spec\BC[X_1,\dots,X_{n-1}]$ to $E_n\simeq
\Spec\BC[\sigma_2,\dots,
\sigma_n]$ defined by $\sigma_i \mapsto (-1)^i\alpha_i$.
\end{defi}

Of course, we have in mind the geometric interpretation of $LL$ explained
below.
Let $v\in V$. The orbit $\overline{v} \in W\cq V$ is represented by a pair
$$(x,y) \in \BC\times Y,$$
where $x=f_n(v)$ and $y$ is the point in $Y$ with coordinates
$(f_1(v),\dots,f_{n-1}(v))$
(this coding of points in $W\cq V$ will be used thoughout this article).
Denote by $p$ the projection $W\cq V \rightarrow Y, (x,y) \mapsto y$.
\begin{defi} For any point $y$ in $Y$, we denote
by $L_y$ the fiber of the projection
$p:W\cq V \rightarrow Y$ over $y$.
\end{defi}
For any $y\in Y$, the affine line $L_y$
intersects the discriminant $\CH$ in $n$ points
(counted with multiplicities), whose coordinates are
$$(x_1,y),\dots,(x_n,y),$$
where $\{\{x_1,\dots,x_n\}\}$ is the multiset of solutions in $X_n$
of $\Delta_f=0$
where each $\alpha_i$ has been replaced by its value at $y$.
We have $$LL(y) = \{\{x_1,\dots,x_n\}\}.$$

\begin{defi}
The \emph{bifurcation locus} is the algebraic hypersurface $\CK\subseteq Y$
defined by the equation $\Disc_{X_n}(\Delta_f)=0$.
\end{defi}

By $\Disc_{X_n}(\Delta_f)$, we mean the discriminant of $\Delta_f$, seen
as a one-variable polynomial with coefficients in $\BC[Y]$.

The first theorem 
of this section is an (elementary) extension to our context of
a result from \cite{looijenga}.

\begin{theo}
\label{theocovering}
The morphism $LL$ is a ramified covering of
degree $n!h^n/|W|$. It restricts to an unramified covering 
$Y-\CK \twoheadrightarrow E_n^{\reg}$.
\end{theo}

The first lemma generalises
\cite[Theorem 1.4]{looijenga}:

\begin{lemma}
$LL$ is \'etale on $Y-\CK$.
\end{lemma}

\begin{proof}
As mentioned in \cite[(1.5)]{looijenga}, the lemma will follow
if we prove that for all $y\in Y-\CK$ with $LL(y)=\{x_1,\dots,x_n\}$,
the hyperplanes $H_1,\dots,H_n$ tangent to $\CH$ at the $n$
distinct points $(x_1,y),\dots,(x_n,y)$ are in general position.

To prove this, we use Proposition \ref{span}.
Each $H_i$ is spanned by $\overline{\xi}_1(x_i,y),\dots,
\overline{\xi}_n(x_i,y)$. Let $(\varepsilon_1,\dots,\varepsilon_n)$
be the basis of $(W \cq V)^*$ dual to
 $(\frac{\partial}{\partial X_1},\dots,
\frac{\partial}{\partial X_n})$. Let $l_i = \sum_j \lambda_{j,i} \varepsilon_j$
be a non-zero vector in $(W \cq V)^*$ orthogonal
to $H_i$. This amounts to taking a non-zero
column vector $(\lambda_{j,i})_{j=1,\dots,n}$
in the kernel of $M(x_i,y)$, or equivalently an eigenvector of $M_0(y)$
associated to the eigenvalue $-x_i$. By assumption, the $x_i$ are distinct.
The eigenvectors are linearly independent.
\end{proof}

The theorem then follows from:

\begin{lemma}
$LL$ is finite, of degree $n!h^n/|W|$.
\end{lemma}

\begin{proof}
Using the first lemma, we observe that $\alpha_2,\dots,\alpha_n$
are algebraically independent.
The lemma follows from a standard argument from the theory of
quasi-homogeneous morphisms. Since each $X_i$ has weight $d_i$,
the weighted Poincar\'e series of
$\BC[Y]=\BC[X_1,\dots,X_{n-1}]$ is $\prod_{i=1}^{n-1} \frac{1}{1-t^{d_i}}$.
The weighted Poincar\'e series of
$\BC[\alpha_2,\dots,\alpha_n]$ is $\prod_{i=2}^n \frac{1}{1-t^{ih}}$.
The extension is finite, of degree 
$$\prod_{i=1}^{n-1} \frac{(i+1)h}{d_i} = \frac{n!h^{n-1}}{d_1\dots d_{n-1}}=
\frac{n!h^n}{|W|}.$$
\end{proof}

\begin{coro}

The space $Y-\CK$ is a $K(\pi,1)$.
\end{coro}

This corollary suffices to obtain one of the results mentionned in the
introduction, Theorem \ref{theointropres}.
Indeed, it was observed in the last paragraph of \cite[section 4]{BM}
that proving that $\pi_2(Y-\CK)=0$ was sufficient to turn our conjectures
into theorems. At this stage, the theorem still relies on the brutal
computations described in \emph{loc. cit.}. A more direct approach
is possible:

Since $\pi_2(Y-\CK)=0$, we may directly 
compute a presentation for $B$
using a Van Kampen type method.
The map $p:W\cq V^{\reg} \rightarrow Y$ restricts to a fibration,
still denoted by $p$, between $E:=p^{-1} (Y-\CK)$ and $Y -\CK$.
Choose a basepoint $y\in Y-\CK$, let $F$ be the fiber $L_y\cap W\cq V$,
choose a basepoint $(x,y)\in F$.
Consider the following commutative diagram, whose first line is the
end of the fibration exact sequence, and where $\beta$ comes from the inclusion
of spaces:
$$\xymatrix{ \pi_2(Y-\CK,y) = 1
\ar[r]  & \pi_1(F,(x,y)) \ar[r]^{\iota} \ar[dr]_{\alpha} 
& \pi_1(E,(x,y)) \ar[r]^{p_*} \ar[d]_{\beta} & \pi_1(Y -\CK,y) \ar[r] & 1 \\
   & & \pi_1(W\cq V^{\reg},(x,y)).}$$
We are in the context of \cite[Theorem 2.5]{zariski}, from which we
conclude that $\alpha$ is surjective.
The space $E$ is the complement in $W\cq V \simeq \BC^n$ of the hypersurface
$\CH \cup p^{-1}(\CK)$, where $\CH$ is the discriminant (defined by
$\Delta(W,f)=0$). The space $W\cq V^{\reg}$ is the complement of $\CH$.
By \cite[Lemma 2.1]{zariski}, the kernel of $\beta$ is 
generated by the generators-of-the-monodromy around $p^{-1}(\CK)$.

Let us prove that the exact sequence in the first line of the diagram is
split.
Let $\gamma_1,\dots,\gamma_m$ be generators-of-the-monodromy in 
$\pi_1(Y -\CK,y)$ such that $\{\gamma_1,\dots,\gamma_m\}$ generates
$\pi_1(Y-\CK,y)$.
Let $\CR$ be a finite set of relations between $\gamma_1,\dots,\gamma_m$
such that
$$\pi_1(Y-\CK,y)) \simeq \left< \gamma_1,\dots,\gamma_m \left|
\CR \right. \right>$$
(fundamental groups of hypersurface complements are always
finitely presentable; a way to prove this is by intersecting with 
a generic $2$-plane and applying Van Kampen method). Choose actual
paths representing the generators and choose homotopies realising the relations.
All these homotopies occur within a compact subspace $C\subseteq Y$.
The image $LL(C)$ is compact. Choose $x$ a complex number with sufficiently
large module, such that the horizontal section $\{x\} \times C$ does
not contain any point of $\CH$. Let $\tilde{\gamma}_1,\dots,
\tilde{\gamma}_m$ be the elements of $\pi_1(E,(x,y))$
 obtained by lifting in $\{x\}\times C$ the paths representing
 ${\gamma_1},\dots, {\gamma_m}$. These elements satisfy the same
defining relations (where each $\gamma_j$ is replaced by $\tilde{\gamma}_j$ --
we denote this set of relations by $\tilde{\CR}$).
The map $\gamma_j \mapsto \tilde{\gamma}_j$ is a section of $p_*$.

The fiber $F$ identifies with $\BC - LL(y)$, that is, since $y\notin\CK$,
the complement of $n$ points in $\BC$. Its fundamental group is a free
group of rank $n$. A standard generating system $f_1,\dots,f_n$ consisting
of generators-of-the-monodromy is obtained when choosing a \emph{planar
spider} (in the sequel, we will use a simpler way of encoding
standard generators, by means of \emph{tunnels}).

The group $\pi_1(E,(x,y))$ is generated by $f_1,\dots,f_n,\tilde{\gamma}_1,
\dots,\tilde{\gamma}_m$. Each $\gamma_j$ defines a monodromy automorphism
$\phi_j \in \Aut( \pi_1(F,(x,y)))$ such that, for all $i\in \{1,\dots,n\}$ and
all $j\in \{1,\dots,m\}$, one has
$$\tilde{\gamma}_j f_i \tilde{\gamma}_j^{-1} = \phi_j(f_i).$$
The semi-direct product structure implies that we have a presentation
$$\pi_1(E,(x,y)) \simeq \left< f_1,\dots,f_n,\tilde{\gamma}_1,
\dots,\tilde{\gamma}_m \left|
\tilde{\gamma}_j f_i \tilde{\gamma}_j^{-1} = \phi_j(f_i),
\tilde{\CR} \right. \right>.$$
Since $\ker \beta$ is the normal subgroup generated the $\tilde{\gamma}_j$'s,
we obtain
$$\pi_1(W\cq V^{\reg},(x,y)) \simeq
\left< f_1,\dots,f_n,\tilde{\gamma}_1,\dots,\tilde{\gamma}_m \left|
\tilde{\gamma}_j f_i \tilde{\gamma}_j^{-1} = \phi_j(f_i),
\tilde{\CR}, \tilde{\gamma}_1=1,\dots,\tilde{\gamma}_m=1 \right. \right>.$$
Since the relations in $\tilde{\CR}$ do not involve $f_1,\dots,f_n$, the
presentation may be rewritten as:

\begin{theo}
\label{theovankampen}
Let $\phi_1,\dots,\phi_m$ be monodromy automorphisms of the fundamental
group of $L_y \cap W \cq V^{\reg}$ associated with
generators $\gamma_1,\dots,\gamma_m$
of the fundamental group of $Y-\CK$.
The braid group of $W$ admits the presentation
$$\pi_1(W\cq V^{\reg},(x,y)) \simeq
\left< f_1,\dots,f_n \left|
f_i = \phi_j(f_i) \right. \right>,$$
where one takes a relation $f_i = \phi_j(f_i)$
for each pair $i,j$ with
$i\in\{1,\dots,n\}$ and $j\in \{1,\dots,m\}$.
\end{theo}

The above theorem refines the main result from \cite{zariski},
which unfortunately was written without knowledge of classical article
by Looijenga \cite{looijenga}.

Later on, we will be able to write presentations of $B(W)$ in much more
explicit ways.

\section{Tunnels, labels and Hurwitz rule}

Let $W$ be an irreducible well-generated complex reflection group.
We keep the notations from the previous section.
Let $y\in Y$.
Let $U_y$ be the complement in $L_y$ of the vertical imaginary half-lines below
the points of $LL(y)$, or in more formal terms:
$$U_y := \{(z,y) \in L_y | \forall x \in LL(y),
\re(z)=\re(x) \Rightarrow \im(z)>\im(x)\}.$$

In this example, the support of $LL(y)$ consists of $4$ points,
and $U_y$ is the complement of three half-lines:
$$\xy
(-5,-10)="11",
(-5,0)="1", (-5,6)="2", (7,-4)="3", (7,-10)="33",
(14,2)="4", (14,-10)="44",
(3,4)*{_{U_y}},
"1"*{\bullet},"2"*{\bullet},"3"*{\bullet},"4"*{\bullet},
"11";"1" **@{-},
"2";"1" **@{-},
"33";"3" **@{-},
"44";"4" **@{-}
\endxy $$
Let $$\CU:= \bigcup_{y\in Y} U_y.$$

\begin{lemma}
The subset $\CU$ is dense in $W\cq V^{\reg}$, open and contractible.
\end{lemma}

\begin{proof}
The first two statements are clear.

Define a continuous function $\beta:Y\rightarrow \BR$ by
 $$\beta(y):=\max\{\im(x)|x \in LL(y)\}-1.$$
Points of $W\cq V$ are represented by pairs $(x,y)\in\BC\times Y$,
or equivalently
by triples $(a,b,y)\in \BR\times\BR\times Y$,
where $a=\re(x)$ and $b=\im(x)$.
For $t\in [0,1]$, define $\phi_t: W\cq V \rightarrow W\cq V$ by
$$\phi_t(a,b,y) := \left\{ 
\begin{matrix} (a,b,y) & \text{if } b\geq \beta(y), \\
(a,b+t(\beta(y)-b),y) & \text{if } b \leq \beta(y). \end{matrix}
\right.$$
Each $\phi_t$ preserves $\CU$ and
the homotopy $\phi$ restricts to a deformation retraction of
$\CU$ to $$\bigcup_{y\in Y} \{(x,y) \in L_y | \im(x) \geq \beta(y)\}.$$
The latter is a locally trivial bundle over the contractible space $Y$,
with contractible fibers (the fibers are half-planes). Thus it is contractible.
\end{proof}

As explained in the topological preliminaries, we may
(and will) use $\CU$ as ``basepoint'' for $W \cq V^{\reg}$.

\begin{defi}
The \emph{braid group of $W$} is $B(W):=\pi_1(W \cq V^{\reg},\CU)$.
\end{defi}

As for other notions actually depending on $W$,
 we usually write $B$ instead of $B(W)$,
since most of the time we implicitly refer to a given $W$.

\begin{defi}
A \emph{semitunnel} is a
triple $T=(y,z,L)\in Y\times \BC \times \BR_{\geq 0}$ such
that $(z,y) \in \CU$ and the affine segment $[(z,y),(z+L,y)]$ lies in 
$W\cq V^{\reg}$.
The path $\gamma_T$ associated with $T$ is the path $t\mapsto (z+tL,y)$. 
The semitunnel $T$ is a \emph{tunnel} if in
addition $(z+L,y) \in \CU$.
\end{defi}

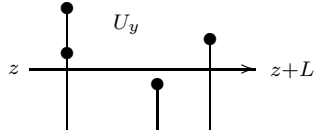
$$\xy
(-5,-10)="11",
(-5,0)="1", (-5,6)="2", (7,-4)="3", (7,-10)="33",
(14,2)="4", (14,-10)="44",
(3,4)*{_{U_y}},
"1"*{\bullet},"2"*{\bullet},"3"*{\bullet},"4"*{\bullet},
"11";"1" **@{-},
"2";"1" **@{-},
"33";"3" **@{-},
"44";"4" **@{-},
(-10,-2);(20,-2) **@{-}, *\dir{>},
(-12,-2)*{_{z}}, (25,-2)*{_{z+L}}
\endxy $$

The distinction between tunnels and semitunnels should be understood
in light of our topological conventions: if $T$ is a tunnel,
$\gamma_T$ represents
an element $$b_T\in \pi_1(W\cq V^{\reg},\CU),$$ 
while semitunnels will be
used to represent points of the universal cover
$(\widetilde{W\cq V^{\reg}})_\CU$ (see Section \ref{sectionuniversalcover}).
 
\begin{defi}
An element $b\in B$ is \emph{simple} if $b=b_T$ for some tunnel $T$.
\end{defi}
 
Each tunnel lives in a given $L_y$, where it may be represented by
an horizontal (= constant imaginery part) segment avoiding $LL(y)$
and with endpoints in $U_y$.
The triple $(y,z,L)$ may be uniquely recovered
from $[(z,y),(z+L,y)]$. A frequent abuse of terminology will
consist of using the term \emph{tunnel} (or \emph{semitunnel})
to designate either the triple $(y,z,L)$, or the segment
$[(z,y),(z+L,y)]$, or the pair $(y,[z,z+L])$, depending on the context
(in particular, when intersecting tunnels with geometric objects,
the tunnels should be understood as affine segments).

Let $y\in Y$.
Let $(x_1,\dots,x_k)$ be the ordered support of $LL(y)$.
The space $p_x((L_y \cap W\cq V^{\reg}) - U_y)$ is an union of
$m$ disjoint open affine intervals $I_1,\dots,I_k$,
where $$I_i := \left\{ \begin{matrix} (x_i-\sqrt{-1}\infty,x_i)
 & \text{if } 
i=1 \text{ or } (i> 1 \text{ and } \re(x_{i-1}) <\re(x_{i})), \\
(x_{i-1},x_i) & \text{otherwise.}
\end{matrix} \right.$$
(by $(x_i-\sqrt{-1}\infty,x_i)$, we mean the open vertical half-line below
$x_i$).
In the first case (when $I_i$ is not bounded),
we say that $x_i$ is \emph{deep}.

$$
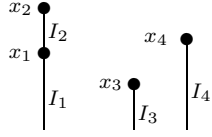
\xy
(-5,-10)="11",
(-5,0)="1", (-5,6)="2", (7,-4)="3", (7,-10)="33",
(14,2)="4", (14,-10)="44",
"1"*{\bullet},"2"*{\bullet},"3"*{\bullet},"4"*{\bullet},
(-8,0)*{_{x_1}},(-8,6)*{_{x_2}},(4,-4)*{_{x_3}},(10,2)*{_{x_4}},
(-3,-6)*{_{I_1}},(-3,3)*{_{I_2}},(9,-8)*{_{I_3}},(16,-5)*{_{I_4}},
"11";"1" **@{-},
"2";"1" **@{-},
"33";"3" **@{-},
"44";"4" **@{-}
\endxy $$

Choose a \emph{system of elementary tunnels} for $y$. By this, we mean the 
choice, for each $i=1,\dots,k$, of a small tunnel $T_i$ in $L_y$
crossing $I_i$
and not crossing the other intervals; let $s_i:=b_{T_i}$ be the associated 
element of $B$.

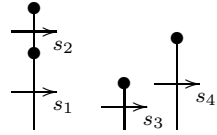
$$\xy
(-5,-10)="11",
(-5,0)="1", (-5,6)="2", (7,-4)="3", (7,-10)="33",
(14,2)="4", (14,-10)="44",
"1"*{\bullet},"2"*{\bullet},"3"*{\bullet},"4"*{\bullet},
"11";"1" **@{-},
"2";"1" **@{-},
"33";"3" **@{-},
"44";"4" **@{-},
(-8,-5);(-2,-5) **@{-}, *\dir{>}, (-1,-7)*{_{s_1}},
(-8,3);(-2,3) **@{-}, *\dir{>}, (-1,1)*{_{s_2}},
(4,-7);(10,-7) **@{-}, *\dir{>}, (11,-9)*{_{s_3}},
(11,-4);(17,-4) **@{-}, *\dir{>}, (18,-6)*{_{s_4}}
\endxy $$

These elements depend only on $y$ and not of the explicit choice
of elementary tunnels.

\begin{defi}
The sequence $lbl(y):=(s_1,\dots,s_k)$ is the \emph{label of $y$}.
Let ${i_1},{i_2},\dots,{i_l}$ be the indices of the successive 
deep points of $LL(y)$.
The \emph{deep label} of $y$ is the subsequence $(s_{i_1},\dots,s_{i_l})$.
\end{defi}

In the above example, the deep label is $(s_1,s_3,s_4)$.
The length of the label is $n$ if and only if $y\in Y-\CK$. The
deep label coincides with the label if and only if $y \in Y^{\gen}$.

It is clear that if $y\in Y^{\gen}$, the elements $s_1,\dots,s_n$ such
that $(s_1,\dots,s_n)=lbl(y)$ are braid reflections.

Later on, it will appear
that the pair $(LL(y),lbl(y))$ uniquely determines $y$
(Theorem \ref{theolabel}). 
It is also possible to characterise, for a given configuration (not
necessarily regular), the possible labels, thus
to explicitly describe of the fibers of $LL$.

Consider the case $y=0$ (given by the equations $X_1=0,\dots,X_{n-1}=0$).
The multiset $LL(y)$ has support $\{0\}$ with multiplicity $n$.

\begin{defi}
We denote by $\delta$ the simple element such that $lbl(0)=(\delta)$.
\end{defi}

This element plays the role Deligne's element $\Delta$.
Let $v\in V^{\reg}$ such that the $W$-orbit $\overline{v}$ lies in 
$L_0$. Brou\'e-Malle-Rouquier consider the element (denoted by $\bm{\pi}$, 
\cite[Notation 2.3]{bmr}) in the pure braid group $P(W):= \pi_1(V^{\reg})
= \ker(B(W) \rightarrow W)$ represented by the loop
$$t \mapsto v\exp(2\sqrt{-1}\pi t).$$
They observed that this element lies in the center of $B$
(\cite[Theorem 2.24]{bmr}).
Since $X_n$ has weight $h$, $\delta^h$ coincides with this element
$\bm{\pi}$. More precisely, $\delta$ is represented by the loop
in $W\cq V^{\reg}$ image of the path in $V^{\reg}$
$$t \mapsto v\exp(2\sqrt{-1}\pi t/h).$$
In particular:

\begin{lemma}
The element $\delta^h$ is central in $B$. The image of
$\delta$ in $W$ is $\zeta_h$-regular, in the sense
of \cite{springer}.
\end{lemma}

Each tunnel lives in a single fiber $L_y$. Simple elements represented by
tunnels in different fibers may be compared since being a tunnel is an
open condition, which is preserved when $y$ is perturbated.

\begin{defi}
Let $T=(y,z,L)$ be a tunnel. A \emph{$T$-neighbourhood} of $y$
is a path-connected neighbourhood $\Omega$ of $y$ in $Y$ such that,
for all $y'\in \Omega$, $T':=(y',z,L)$ is a tunnel.
\end{defi}

Such neighbourhoods clearly exist for all $y\in Y$.

\begin{lemma}[Hurwitz rule]
Let $T=(y,z,L)$ be a tunnel, representing a simple element $s$.
Let $\Omega$ be a $T$-neighbourhood of $y$.
For all $y'\in \Omega$, $T':=(y',z,L)$ represents $s$.
\end{lemma}

\begin{proof}
This simply expresses that the tunnels $(y',z,L)$ and
$(y,z,L)$ represent homotopic paths, which is clear by definition of
$\Omega$.
\end{proof}

\begin{remark}
Let $T_1,\dots,T_k$
is a system of elementary tunnels for $y$. Let $\Omega_i$ be
a $T_i$-neighbourhood for $y$. A \emph{standard neighbourhood} of $y$
could be defined as a path-connected neighbourhood $\Omega$
of $y$ inside
$\cap_{i=1}^k \Omega_i$. These standard neighbourhoods form a basis
for the topology of $Y$.
A consequence of Hurwitz rule is that the label of $y$ may uniquely
be recovered once we know the label of a single $y'\in \Omega$. The converse
is not true, unless $y\in Y -\CK$. Let $y'\in \Omega$, let $S$ and $S'$ 
be the strata of $LL(y)$ and $LL(y')$ (with respect to the discriminant
stratification of $E_n$).
We have $S' \leq S$. The label of $y'$
is uniquely determined by the label of $y$ if and only if $S=S'$.
Otherwise, the ambiguity is caused by the extra ramification at $y$.
\end{remark}

The remainder of this section consists of various consequences of
Hurwitz rule.

\begin{coro}
\label{corodeeplabel}
Let $y\in Y$.
Let $(s_{i_1},\dots,s_{i_l})$ be the deep label of $y$.
We have $s_{i_1}\dots s_{i_l}=\delta$.
\end{coro}

\begin{proof}
Any tunnel $T$ 
deep enough and long enough represents $s_{i_1}\dots s_{i_l}$.
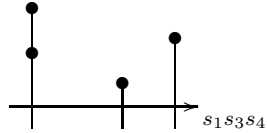
$$\xy
(-5,-10)="11",
(-5,0)="1", (-5,6)="2", (7,-4)="3", (7,-10)="33",
(14,2)="4", (14,-10)="44",
"1"*{\bullet},"2"*{\bullet},"3"*{\bullet},"4"*{\bullet},
"11";"1" **@{-},
"2";"1" **@{-},
"33";"3" **@{-},
"44";"4" **@{-},
(-8,-7);(17,-7) **@{-}, *\dir{>}, (22,-9)*{_{s_1s_3s_4}}
\endxy $$
The origin $0\in Y$ lies in
a $T$-neighbourhood of $y$. To conclude, apply Hurwitz rule.
\end{proof}

\begin{coro}
Let $x \in E_n^{\gen}$.
Let $\beta\in \pi_1(E_n^{\reg},x)$, $y\in LL^{-1}(x)$ and $y':= y \cdot
\beta$. Let $(b_1,\dots,b_n)$ be the label of $y$ and 
$(b_1',\dots,b_n')$ be the label of $y'$.
Then $$(b_1',\dots,b_n') = (b_1,\dots,b_n) \cdot \beta,$$
where $\beta$ acts by right Hurwitz action.
\end{coro}

\begin{proof}
It is enough to prove this for a standard generator $\bm{\sigma_i}$.
Let $(x_1,\dots,x_n)$ be the ordered support of $x$.
By Hurwitz rule, we may adjust the imaginary parts of the $x_i$'s without
affecting the label; in particular we may assume that $\im(x_i) < \im(x_{i+1})$.
We may find tunnels $T_-=(y,z_-,L)$ and $T_+=(y,z_+,L)$ as in the picture
below:
$$
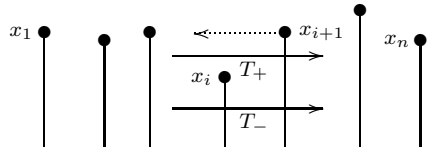
\xy
(-10,0)="1", (-2,-1)="2", (4,0)="3", 
(14,-6)="4", (22,0)="5", (32,3)="6", (40,-1)="7",
(-10,-15)="11", (-2,-15)="22", (4,-15)="33", 
(14,-15)="44", (22,-15)="55", (32,-15)="66", (40,-15)="77",
"1";"11" **@{-}, "2";"22" **@{-}, "3";"33" **@{-}, "4";"44" **@{-},
"5";"55" **@{-}, "6";"66" **@{-}, "7";"77" **@{-},
"1"*{\bullet},"2"*{\bullet},"3"*{\bullet},"4"*{\bullet},
"5"*{\bullet},"6"*{\bullet},"7"*{\bullet},
(11,-6)*{_{x_i}},(27,0)*{_{x_{i+1}}},
(-13,0)*{_{x_1}},(37,-1)*{_{x_{n}}},
"5";(10,0) **@{.}, *\dir{>},
(7,-3);(27,-3) **@{-}, *\dir{>},
(7,-10);(27,-10) **@{-}, *\dir{>},
(18,-12)*{_{T_-}},
(18,-5)*{_{T_+}}
\endxy $$
The path in $E_n^{\reg}$ where $x_{i+1}$ moves along the dotted arrow
and all other points are fixed represents $\bm{\sigma_i}$.
Applying Hurwitz rule to $T_+$, we obtain $b'_i=b_{i+1}$; applying
Hurwitz rule to $T_-$, we obtain $b'_ib'_{i+1}= b_ib_{i+1}$. The result
follows.
\end{proof}

\begin{coro}
\label{corocardhurwitz}
Let $y\in Y^{\gen}$. 
The cardinal of the Hurwitz orbit $lbl(y) \cdot B_n$ is at most $n!h^n/|W|$,
and there is an equivalence between:
\begin{itemize}
\item[(i)] $|lbl(y)\cdot B_n| = n!h^n/|W|$.
\item[(ii)] The map
$Y^{\gen} \rightarrow E_n\times B^n, y \mapsto (LL(y),lbl(y))$
is injective.
\end{itemize} 
\end{coro}

In the next section, we will prove that both (i) and (ii) hold.
This is not a trivial statement.

\begin{proof}
By the previous corollary, one has $lbl(y\cdot\beta)=lbl(y)\cdot\beta$.
Since $y\in Y-\CK$, we have $|y\cdot\beta| = n!h^n/|W|$ (Theorem
\ref{theocovering}). The corollary follows.
\end{proof}

\begin{coro}
\label{coros1si}
Let $s$ be a simple element. There exists $y\in Y^{\gen}$
and $i\in\{1,\dots,n\}$ such that
$s=s_1\dots s_i$,  where
$(s_1,\dots,s_n):=lbl(y)$.
\end{coro}

\begin{proof}
Let $T$ be a tunnel representing $s$.
Any $T$-neighbourhood of $y$ contains generic points.
We may assume that $y\in Y^{\gen}$. The picture below
explains, on an example, how
to move certain points (following the dotted paths) of
the underlying configuration to reach a suitable $y'$:
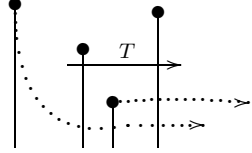
$$\xy
(-5,-13)="11",
(-5,0)="1", (-5,6)="2", (4,0)="3", (4,-13)="33",
(14,5)="4", (14,-13)="44", (8,-7)="5", (8,-13)="55",
"5"*{\bullet},"2"*{\bullet},"3"*{\bullet},"4"*{\bullet},
"33";"3" **@{-},
"11";"2" **@{-},
"55";"5" **@{-},
"44";"4" **@{-},
(2,-2);(17,-2) **@{-}, *\dir{>}, (10,0)*{_{T}},
"2";(10,-10) **\crv{~*=<4pt>{.} (-5,0),(-5,-13)},
(10,-10);(20,-10) **\crv{~*=<4pt>{.} (13,-10),(17,-10)},
*\dir{>},
"5";(26,-7) **\crv{~*=<4pt>{.} (16,-6)},
*\dir{>}
\endxy $$
This path in $E_n^{\reg}$
lifts, via $LL$, to a path in a $T$-neighbourhood of $y$
whose final point $y'$ satisfies the conditions of the lemma.
\end{proof}

\section{Reduced decompositions of Coxeter elements}

Most results on this section follow from the generalisation to well-generated
complex $2$-reflection groups of a property of finite Coxeter groups
initially conjectured by Looijenga, \cite[(3.5)]{looijenga}, and
solved in a letter from Deligne to Looijenga
(citing discussions with Tits and Zagier),
\cite{deligneletter}.
This property is more or less equivalent to
 \cite[Fact 2.2.4]{dualmonoid}, used in the our earlier construction
of the dual braid monoid. It relies on a case-by-case argument, which
we hope to improve in the future.

As in the previous sections, $W$ is a well-generated complex $2$-reflection
group, but not necessarily irreducible.
We study the images of simple elements under
the natural projection $\pi: B \rightarrow W$.
Recall that the projection $\pi$ is part of the fibration exact sequence
$$\xymatrix@1{ 1 \ar[r] & \pi_1(V^{\reg})
\ar[r] & \pi_1(W\cq V^{\reg} ) \ar[r]^{\phantom{mmm}\pi} & W \ar[r] & 1.}$$
For this exact sequence to be be well-defined, one
has to make consistent choices of basepoints in $V^{\reg}$
and in $W\cq V^{\reg}$.
We have already described our ``basepoint''
$\CU$ in $W\cq V^{\reg}$. Choose $u\in \CU$ and choose a 
preimage $\tilde{u}$ of $u$ in $V^{\reg}$. If $u'\in \CU$ is another choice and
if $\gamma$ is a path in $\CU$ from $u$ to $u'$, $\gamma$ lifts
to a unique path $\tilde{\gamma}$ starting at $\tilde{u}$; since
$\CU$ is contractible, the fixed-endpoint 
homotopy class of $\tilde{\gamma}$ (and, in particular, its final point)
does not depend on $\gamma$. In other words, once we have chosen
a preimage of one point of $\CU$, we have a natural section $\tilde{\CU}$
of $\CU$ in $V^{\reg}$, as well as a transitive system of isomorphisms
between $(\pi_1(V^{\reg},\tilde{u}))_{\tilde{u}\in\tilde{\CU}}$.
From now on, we assume we have made
such a choice, and we define the pure braid group
as $\pi_1(V^{\reg},\tilde{\CU})$.
(The $|W|$ possible choices yield conjugate morphisms $\pi$).

Recall that an element $w \in W$ is \emph{$\zeta$-regular}, in the
sense of \cite{springer}, if at admits in $\zeta$-eigenvector in $V^{\reg}$.
We then say that $\zeta$ is a \emph{regular eigenvalue}, and its order $d$ is
a \emph{regular number}.
Springer proved that, if $\zeta$ is a regular eigenvalue, all $\zeta$-regular
elements of $W$ are conjugate in $W$.

\begin{defi}
Let $W$ be an irreducible well-generated complex $2$-reflection group.
An element $c\in W$ is a
\emph{(generalised) Coxeter element} if it is $\zeta_h$-regular.
More generally, if $W$ is a well-generated complex $2$-reflection group
decomposed as a sum $W=\bigoplus_i W_i$ of irreducible groups,
a Coxeter element in $W$ is a product
$c=\prod_i c_i$ of Coxeter elements in each $W_i$.
\end{defi}

A consequence of our description of $\delta$ (see previous section) is
that $c:=\pi(\delta)$ is a Coxeter element.
The other Coxeter elements, which are conjugate to $c$,
appear when considering other basepoints over $\CU$.

We denote by $R$ the set of all reflections in $W$.
As in \cite{dualmonoid}, for all $w\in W$, we denote by $\Red_T(w)$ the
set of \emph{reduced $R$-decompositions of $w$}, that it is minimal
length sequences of elements of $R$ with product $w$.
Since $R$ is closed under conjugacy, $\Red_T(c)$ is stable under Hurwitz
action.
We also consider the length function $l_R:W\rightarrow \BZ_{\geq 0}$,
whose value at $w$ is the common length of the elements of $\Red_R(w)$.
We also consider the two partial orderings of $W$
defined as follows: for all $w,w' \in W$, we set
$$w \preccurlyeq_R w' :\Leftrightarrow l_R(w) + l_R(w^{-1}w')  = l_R(w)$$
and 
$$w' \succcurlyeq_R w :\Leftrightarrow l_R(w') + l_R(w'w^{-1})  = l_R(w).$$ 
Since $R$ is invariant by conjugacy, we have
$w\preccurlyeq_R w' \Leftrightarrow w' \succcurlyeq_R w$.

It is worth noting that the basic Lemma 1.2.1. of \cite{dualmonoid} does not
hold when $W$ is not a complexified real group: for obvious reasons,
one has $l_R(w)\geq n$ for all $w\in W$, but sometimes
$l_R(w) > n$.
However, we have:

\begin{lemma}
\label{lemmaparacox}
Let $c$ be a Coxeter element in a well-generated complex $2$-reflection group
$W$. Let $w\in W$ with $w\preccurlyeq_R c$.
Then
\begin{itemize}
\item $\codim \ker(w-1)= l_R(w)$,
\item $w$ is a Coxeter
element in the group $W_w:=\{w'\in W | \ker(w-1) \subseteq \ker(w'-1) \}$,
\item $\Red_{R}(w)= \Red_{R\cap W_w}(w)$.
\end{itemize}
\end{lemma}

That $W_w$ is a reflection group is classical (the fixator of any subspace
of $V$ is a reflection group, called \emph{parabolic subgroup}).
One checks (e.g., on the classification) that parabolic subgroups of 
well-generated complex reflection groups are well-generated.

\begin{proof}
The lemma reduces to the case when $W$ is irreducible, which we assume
from now on.

(i) Assume that $w=c$. We have $\ker(c-1)=0$ since the eigenvalues of $c$
are $\zeta_h^{1-d_1},\dots,\zeta_h^{1-d_n}$ and $1< d_i \leq h$ for all $i$.
That $l_R(c)=n$ is easy and classical (it
follows, for example, from the main theorem in \cite{zariski} --
in this reference, the statements
 about the eigenvalues contain a typo: $\zeta$ should be replaced by
$\overline{\zeta}$).

(ii) Let $w\prec c$. 
Let us check that $w$ is a Coxeter element in $W_w$.
When $W$ is complexified real, this follows from
\cite[Lemma 1.4.3]{dualmonoid}.
For $G(e,e,r)$, it follows from \cite[Proposition 4.3]{BC}.
The remaining exceptional cases
are easy (an induction reduces the problem to the case $l_R(w)=n-1$) to
check (though infortunately we don't have a general argument).

(iii) We have $w=c_1\dots c_k$, where the $c_i$ are Coxeter elements of 
irreducible parabolic subgroups. By (i), we have, for each $c_i$,
$l_{R_{c_i}}(c_i)= \codim\ker(c_i-1)$. We deduce that
$l_R(w) \leq \sum \codim\ker(c_i-1) = \codim\ker (w-1)$. The opposed 
inequality being obvious, we have $l_R(w) = \codim\ker (w-1)$.

(iv) If $w=r_1\dots r_k$ with $r_i\in R$ and $k=\codim\ker(w-1)$, one
necessarily has $\ker(w-1) = \bigcap_{i=1}^k \ker (r_i-1)$.
In other words, $r_i\in W_w$. This proves the last statement.
\end{proof}

\begin{prop}
\label{propcardhurwitz}
Let $W$ be a well-generated complex $2$-reflection group.
Let $c$ be a generalised Coxeter element in $W$.
Let $w\in W$ with $w\preccurlyeq c$.
The set $\Red_R(w)$ is a
single Hurwitz orbit.
If $W$ is irreducible with Coxeter number $h$,
we have $|\Red_R(c)| = n!h^n/|W|$.
\end{prop}

\begin{proof}
Using Lemma \ref{lemmaparacox},
the proposition reduces to the case when $W$ is irreducible
and $w=c$: in the reducible case,
reduced decompositions of Coxeter elements are ``shuffles'' of reduced
decompositions of the Coxeter elements of the irreducible summands.

The complexified real case is studied in \cite{deligneletter}
(transitivity is easy and does not require case-by-case, see for
example \cite{dualmonoid}, Proposition 1.6.1).
The $G(e,e,r)$ case is a combination of Proposition 6.1 (transitivity)
and Theorem 8.1 (cardinality) from \cite{BC}. The remaining exceptional
types are done by computer.
\end{proof}

We now assume that $W$ is an irreducible well-generated complex $2$-reflection
group and return to the situation and notations of the previous sections.

\begin{theo}
\label{theolooijenga}
Let $y\in Y^{\gen}$, with label $(s_1,\dots,s_n)$.
For $i=1,\dots,n$, let $r_i:=\pi(s_i)$.
The map $\pi^n$ induces an isomorphism of $B_n$-sets between the
Hurwitz orbits of $(s_1,\dots,s_n)$ and of $(r_1,\dots,r_n)$.
\end{theo}

\begin{proof}
We have $\delta = s_1\dots s_n$. Let $c:=r_1\dots r_n$.
As mentioned above, $c$ is a Coxeter element (\cite[Lemma 3.1]{zariski}).
By Proposition \ref{propcardhurwitz}, the cardinal of the Hurwitz
$(r_1,\dots,r_n)\cdot B_n$ is $n!h^n/|W|$.
By Corollary \ref{corocardhurwitz}, the cardinal of the Hurwitz
orbit $(s_1,\dots,s_n)\cdot B_n$ is at most $n!h^n/|W|$.
One concludes that $|(s_1,\dots,s_n)\cdot B_n| = n!h^n/|W|$ and that
$\pi^n:(s_1,\dots,s_n)\cdot 
B_n \rightarrow (r_1,\dots,r_n)\cdot B_n$ is an isomorphism.
\end{proof}

\begin{coro}
\label{corosimples}
The restriction of $\pi$ to the set $S$ of simple elements is
injective.
\end{coro}

\begin{proof}
Let $s$ and $s'$ be simples such that $\pi(s)=\pi(s')$.
By Corollary \ref{coros1si}, we may find $y,y'\in Y^{\gen}$,
with $lbl(y)=(s_1,\dots,s_n)$ and 
$lbl(y')=(s_1',\dots,s_n')$,
and $i\in\{0,\dots,n\}$ such that $s=s_1\dots s_i$
and $s'=s'_1\dots s'_i$. 
By Proposition \ref{propcardhurwitz} applied to $w=\pi(s)$,
we conclude that we may find $\beta\in B_i$ such that
 $(\pi(s_1),\dots,\pi(s_i))\cdot \beta = (\pi(s'_1),\dots,\pi(s'_i))$.
Similarly, we find $\beta'\in B_{n-i}$ such that
 $(\pi(s_{i+1}),\dots,\pi(s_n))\cdot \beta' = (\pi(s'_{i+1}),\dots,
\pi(s'_n))$. View $B_i \times B_{n-i}$ as a subgroup of $B_n$ (the
first factor braids the first $i$ strings, the second factors braids
the $n-i$ last strings), and set $\beta'':=(\beta,\beta')\in B_n$.
We have $\pi^n(lbl(y))\cdot \beta''=\pi^n(lbl(y'))$.
Applying the theorem, this implies that $lbl(y) \cdot \beta'' = lbl(y')$.
Clearly, $\beta''$ does not modify the product of
the first $i$ terms of the labels.
Thus $s=s'$.
\end{proof}

In the following theorem, $B^*$ denotes the set of finite sequences of elements
of $B$.

\begin{theo}
\label{theolabel}
The map $Y\rightarrow E_n \times B^*, y\mapsto (LL(y),lbl(y))$ is injective.
\end{theo}

The image of the map may easily be described in terms of simple axioms
relating $lbl(y)$ to $LL(y)$. This description is not needed here.

\begin{proof}
The injectivity on $Y^{\gen}$ follows from Corollary \ref{corocardhurwitz}
and the proof of the previous theorem.

Let $y,y'\in Y$ with $LL(y)=LL(y')$ and $lbl(y)=lbl(y')$.
We have to prove that $y=y'$.
Let $(x_1,\dots,x_k)$ and $(n_1,\dots,n_k)$ be the
ordered support and multiplicities of $LL(y)$.
Let $(t_1,\dots,t_k):=lbl(y)$.

For simplicity, we begin with the case when $\{x_1,\dots,x_k\}\in E_k^{\gen}$,
\ie, we assume that $\re(x_1) < \re(x_2) < \dots < \re(x_k)$.
Let $(y_m)_{m\geq 0}$ be 
sequences of points in $Y^{\gen}$ converging
to $y$. 
Let $(s_1^{(m)},\dots,s_n^{(m)}):=lbl(y_m)$.
For $m$ large enough, Hurwitz rule implies that 
$s_1^{(m)}\dots s_{n_1}^{(m)} = t_1$, 
$s_{n_1+1}^{(m)}\dots s_{n_1+n_2}^{(m)}=t_2$, etc...
Similarly, choose a sequence $(y'_m)$ of points in $Y^{\gen}$
converging to $y'$;
one has the corresponding properties for
$(s_1^{(m)'},\dots,s_n^{(m)'}):=lbl(y_m)$.
For $m$ large, $(s_1^{(m)},\dots,s_{n_1}^{(m)})$
and $(s_1^{(m)'},\dots,s_{n_1}^{(m)'})$ are two reduced $R$-decompositions
of the simple $t_1$. By Proposition \ref{propcardhurwitz}, they
are in the same $B_{n_1}$ Hurwitz orbit.
Similarly, $(s_{n_1+1}^{(m)},\dots,s_{n_1+n_2}^{(m)})$
and $(s_{n_1+1}^{(m)'},\dots,s_{n_1+n_2}^{(m)'})$ are in the same
$B_{n_2}$ Hurwitz orbit, etc... We have proved that, for $m$ large enough,
$lbl(y_m)$ and $lbl(y'_m)$ are Hurwitz transformed by a braid
$\beta_m$ which lies in the natural subgroup $B_{n_1}\times B_{n_2} \times \dots
\times B_{n_k}$ of $B_n$. One concludes noting that this subgroup belongs to
the ramification group for the $B_n$-action on $y$.

The case when $\{x_1,\dots,x_k\}\notin E_k^{\gen}$ may be reduced
to the previous case by some suitable perturbation.
\end{proof}

We will also need the following lemma:

\begin{lemma}
\label{lemmaunique}
Let $y\in Y$. 
Let $T=(y,z,L)$ and $T'=(y,z',L')$ be two tunnels in $L_y$
such that $b_{T}=b_{T'}$.
Then $T$ and $T'$ cross the same intervals among $I_1,\dots,I_n$
(in other words, $T$ and $T'$ are homotopic as tunnels drawn in 
$L_y$).
\end{lemma}

\begin{proof}
Up to perturbating $y$, we may assume that $y\in Y^{\gen}$.
Let $(s_1,\dots,s_n):=lbl(y)$, and $(r_1,\dots,r_n):=\pi^n(lbl(y))$.
Let $i_1,\dots,i_l$ (resp. $j_1,\dots,j_m$) be the successive indices
of the intervals among $I_1,\dots,I_n$ crossed
by $T$ (resp. $T'$). We have $b_T=s_{i_1}\dots s_{i_l}$ and $b_T=
s_{j_1}\dots s_{j_m}$.
Assuming that $b_T=b_{T'}$, we obtain
$s_{i_1}\dots s_{i_l} = s_{j_1}\dots s_{j_m}$
and $r_{i_1}\dots r_{i_l} = r_{j_1}\dots r_{j_m}$.
Let $w:=r_{i_1}\dots r_{i_l}$.
By Lemma \ref{lemmaparacox}, $l=m$ and we have 
$$(*) \qquad \ker(w-1) = \bigcap_{k=1}^l \ker(r_{i_k}-1) =
 \bigcap_{k=1}^l \ker(r_{j_k}-1).$$
Assume that $(i_1,\dots,i_l) \neq (j_1,\dots,j_l)$. We may find
$j\in \{j_1,\dots,j_l\}$ such that, for example,
$i_1 < \dots < i_k < j < i_{k+1} < \dots i_l$.
Noting that the element $r_{i_1}\dots r_{i_k} r_jr_{i_{k+1}} \dots r_{i_l}$
is a parabolic Coxeter element, we deduce that $\ker(r_j-1) \not \supseteq
\ker(w-1)$. This contradicts $(*)$.
\end{proof}

\section{The dual braid monoid}

Here again, $W$ is an irreducible well-generated complex $2$-reflection
groups. We keep the notations from the previous sections.
Recall that an element $b\in B$ is \emph{simple} if $b=b_T$ for some tunnel
$T$, and that the set of simple elements is denoted by $S$.

\begin{defi}
The \emph{dual braid monoid} is the submonoid $M$ of $B$ generated by $S$.
\end{defi}

Consider the binary relation $\preccurlyeq$ defined on $S$ as follows.
Let $s$ and $s'$ be simple elements. We write $s \preccurlyeq s'$ if and
only if there exists $(z,y)\in W\cq V^{\reg}$, $L,L'\in \BR_{\geq 0}$ with
$L \leq L'$ such that $(y,z,L)$ is a tunnel representing $s$ and
$(y,z,L')$ is a tunnel representing $s'$.

This section is devoted to the proof of:

\begin{theo}
\label{theogarside}
The monoid $M$ is a Garside monoid, with set of simples $S$ and Garside
element $\delta$. The relation $\preccurlyeq$ defined above on $S$ is the
restriction to $S$ of the left divisibility order in $M$.

The monoid $M$ generates $B$, which inherits a structure of Garside group.
\end{theo}

Several results mentioned in the introduction follow from this theorem:

Theorem \ref{theointrogarside} does not assume irreducibility,
but follows immediately from the irreducible case since
direct products of Garside groups are Garside groups.

Theorem \ref{theocenter} also follows, with arguments similar to
\cite[p. 301]{deligne}: the center of $B$ is generated by
$\delta^a$, where $a$ is the smallest positive integer such that
$c^a$ is central in $W$ (using standard techniques, 
one may also prove that all
$h$-th roots of $\bm{\pi}$ in $B$ are conjugate to $\delta$, and
compute the centralisers of the powers of $\delta$).

Another important consequence of Theorem \ref{theogarside} is that
one may associate a practical contractible simplicial complex acted on
by $B$.
To explain this, we need some additional terminology, imitating
\cite{bestvina}
and \cite{cmw}:

\begin{defi} Let $G$ be a group, let $\Sigma$ be a subset of $G$. We
denote by $\Flag(G,\Sigma)$ the simplicial complex as follows:
\begin{itemize}
\item[($0$)] the set of vertices is $G$,
\item[($1$)] the $1$-simplices are the pairs $\{g,h\}$ of distinct
elements of $G$ such that $g^{-1}h \in \Sigma$ or $h^{-1}g\in \Sigma$,
\item[($k$)] for all $k>1$,
a $(k+1)$-tuple $\{g_0,\dots,g_k\}$ of distinct elements of $G$
is a $k$-simplex if and only if $\{g_i,g_j\}$ is a $1$-simplex
for all $i\neq j$.
\end{itemize}
\end{defi}

In other words, $\Flag(G,\Sigma)$ is the only flag complex
whose $1$-skeleton is the Cayley graph of $(G,\Sigma)$.

The following result is due, in this formulation,
to Charney-Meier-Whittlesey. It is inspired by the work of Bestvina,
\cite{bestvina}, on
spherical type Artin monoid (which itself is based on the work of Deligne).
It was also implicitly used in several articles by T. Brady (e.g.
\cite{brady}).

\begin{theo}[\cite{cmw}]
Let $G$ be a Garside group with set of simples $\Sigma$.
Then $\Flag(G,\Sigma)$ is contractible.
\end{theo}

\begin{coro}
\label{coroflag}
The simplicial complex $\Flag(B,S)$ is contractible.
\end{coro}

The strategy of proof of Theorem \ref{theogarside} is
very similar to the one in \cite{BC}.

\begin{lemma}
Let $s \in S$. We have $l(s)=l_R(\pi(s))$.
\end{lemma}

\begin{prop}
\label{propiiiii}
Let $s,s'\in S$. The following statements are equivalent:
\begin{itemize}
\item[(i)] $s \preccurlyeq s'$,
\item[(ii)] $\exists s''\in S, ss''=s'$,
\item[(iii)] $\pi(s) \preccurlyeq_R \pi(s')$.
\end{itemize}
\end{prop}

\begin{proof}
(i) $\Rightarrow$ (ii). 
Assume that $s \preccurlyeq s'$ and choose $(x,y)\in W\cq V^{\reg}$ 
and $L \leq L'$ such that $(y,x,L)$ is a tunnel representing $s$ and
$(y,x,L')$ is a tunnel representing $s'$.
Then $(y,x+L,L'-L)$ is a tunnel representing $s''\in S$ such
that $ss''=s'$.

(ii) $\Rightarrow$ (iii). The natural length function $l$ is additive
on $B$. Thus, under (ii), we have
 $l(s)+l(s'')=l(s')$. On the other hand, for all $\sigma\in S$,
$l(\sigma)=l_R(\pi(\sigma))$. Thus $l_R(\pi(s))+l_R(\pi(s''))=l_R(\pi(s'))$.
Since $\pi(s'') = \pi(s)^{-1}\pi(s')$, this implies that $\pi(s)\preccurlyeq_R
\pi(s')$.

(iii) $\Rightarrow$ (i).
We may find $y,y'\in Y^{\gen}$, with $lbl(y)=(s_1,\dots, s_n)$ and
$lbl(y')=(s'_1, \dots,s'_n)$ such that
$s=s_1s_2\dots s_{l(s)}$ and $s'=s'_1s'_2\dots s'_{l(s')}$.
Set $w:=\pi(s)$, $w':=\pi(s')$ and, for $i=1,\dots,n$,
set $r_i:=\pi(s_i)$ and $r'_i:=\pi(s'_i)$.
Assuming (iii), we may find $r''_{1},\dots,r''_{l(s')-l(s)}\in R$ such
that $r_1 \dots r_{l(s)} r''_{1}\dots r''_{l(s')-l(s)} = 
r'_1\dots r'_{l(s')}$.
The sequences $$(r_1, \dots, r_{l(s)}, r''_{1},\dots, r''_{l(s')-l(s)})$$
and $$(r'_1,\dots ,r'_{l(s')})$$ both lie in $\Red_R(w')$.
Since $w'\preccurlyeq c$, both sequences lie in the same Hurwitz orbit
(Proposition \ref{propcardhurwitz}). Thus
$$(r'_1,\dots ,r'_{n})$$ and
$$(r_1, \dots, r_{l(s)}, r''_{1},\dots,
 r''_{l(s')-l(s)},r'_{l(s')+1},\dots,
r'_n)$$ are transformed one onto the other by Hurwitz action of
a braid $\beta\in B_n$ only braiding the first $l(s')$ strands.
The Hurwitz tranformed of
$$lbl(y') = (s'_1,\dots,s'_n)$$
by $\beta$ is the label 
$$(s''_1,\dots,s''_n)$$ of some $y''\in Y^{\gen}$.
Since the braid only involves the first $l(s')$ strands,
$s''_1\dots s''_{l(s')}= s'_1\dots s'_{l(s')} = s'$.
One has $\pi(s''_i)= r_i$ for $i=1,\dots,l(s)$,
thus $\pi(s''_1\dots s''_{l(s)}) = \pi(s_1\dots s_{l(s)})$.
By Corollary \ref{corosimples}, this implies
$s''_1\dots s''_{l(s)} = s_1\dots s_{l(s)} =s$.
For $x$ with small enough real and imaginary parts, one may find real
numbers $L,L'$ with $0< L \leq L'$ such that $(y'',x,L)$ is a tunnel with
associated simple $s''_1\dots s''_{l(s)}=s$ and $(y'',x,L')$ is a tunnel with
associated simple $s''_1\dots s''_{l(s')} =s'$.
\end{proof}

\begin{prop}
\label{posetiso}
The map $\pi$ restricts to an isomorphism
$(S,\preccurlyeq)\stackrel{\sim}{\rightarrow}
([1,c],\preccurlyeq_R)$. In particular, $\preccurlyeq$ is an order
relation on $S$.
\end{prop}
 
\begin{proof}
The previous proposition, applied to $s'=\delta$, proves that
$\pi(S)\subseteq [1,c]$. It also proves that $\pi$ induces a morphism
of sets with binary relations $(S,\preccurlyeq){\rightarrow}
([1,c],\preccurlyeq_R)$.
The injectivity is Corollary \ref{corosimples}.

Surjectivity: Choose $y\in Y^{\gen}$.
Let $(s_1,\dots,s_n):=lbl(y)$. Let $r_i:=\pi(s_i)$.
We have $(r_1,\dots,r_n) \in \Red_R(c)$
Let $w \in [1,c]$. We may find $(r'_1,\dots,r'_n) \in \Red_R(c)$ such
that $r'_1\dots r'_{l_T(w)} = w$.
By Proposition \ref{propcardhurwitz}, $(r'_1,\dots,r'_n)$ is Hurwitz
transformed of $(r_1,\dots,r_n)$, thus there exists $y'\in Y^{\gen}$
such that $\pi_*(lbl(y'))=(r'_1,\dots,r'_n)$. The simple element which
is the product of the first $l(w)$ terms of $lbl(y')$ is in $\pi^{-1}(w)$.
\end{proof}

We have the key lemma:

\begin{lemma}
\label{lemmalattice}
The poset $([1,c],\preccurlyeq_R)$ is a lattice.
\end{lemma}

\begin{proof}
The real case is done in \cite{dualmonoid}. The $G(e,e,r)$ is done
in \cite{BC}. The remaining cases may be done by computer.
\end{proof}

\begin{defi}
The \emph{dual braid relations} with respect to $W$ and $c$ are all the
formal relations of the form $$rr'=r'r'',$$
where $r,r',r''\in R$ are such that $r\neq r'$, $rr'\in [1,c]$, and the
relation $rr'=r'r''$ holds in $W$.
\end{defi}

Clearly, dual braid relations only involve reflections in $R\cap [1,c]$.
When $W$ is complexified real, $R\subseteq [1,c]$. This does not hold
in general (see the tables at the end of the article).

\begin{lemma}
\label{lemmadualbraidrelations}
Let $R_c:=R\cap [1,c]$.
The monoid $\mathbf{M}(P_c)$ admits the monoid presentation 
$$\mathbf{M}(P_c) \simeq \left< R_c \left|
\text{dual braid relations} \right. \right>.$$
\end{lemma}

\begin{remark}
When viewed as a group presentation, the presentation of the lemma is 
a presentation for $\mathbf{G}(P_c)$. As soon as we prove
 Theorem \ref{theogarside}, Lemma \ref{lemmadualbraidrelations}
will give an explicit presentation for $B$. A way to reprove
Theorem \ref{theointropres} is by simplifying the (redundent) presentation
given by the lemma. This does not involve any monodromy computation.
\end{remark}

\begin{proof}[Proof of Lemma \ref{lemmadualbraidrelations}]
By definition, $\mathbf{M}(P_c)$ admits the presentation with
generators $R_c$ and a relation
$r_1\dots r_k= r'_1\dots r'_k$ for each pair $(r_1,\dots,r_k)$,
$(r'_1,\dots,r'_k)$ of reduced $R$-decompositions of the same element
$w\in [1,c]$.
Call these relations \emph{Hurwitz relations}.
By transitivity of the Hurwitz action on $\Red_R(c)$ (Proposition
\ref{propcardhurwitz}), the Hurwitz relations are consequences of the 
dual braid relations.
The dual braid relations clearly hold in $\mathbf{M}(P_c)$
(to see this, complete $(r,r')$ to an element
$(r,r',r_3,\dots,r_n)\in\Red_R(c)$). This proves the lemma.
\end{proof}

\begin{proof}[Proof of Theorem \ref{theogarside}]
Set $R_c := R \cap [1,c]$.
We are in the situation of subsection 0.4 in \cite{dualmonoid},
$(W,R_c)$ is a \emph{generated group} and $c$ is \emph{balanced}
(one first observes that $\{w | w \preccurlyeq_{R_c} c \} = 
\{w | w \preccurlyeq_{R} c \}$ and
$\{w | c \succcurlyeq_{R_c} w \} = 
\{w | c \succcurlyeq_{R} w \}$; one concludes noting that $c$ is balanced
with respect to $(W,R)$, which is immediate since $R$ is invariant by
conjugacy).
By Lemma \ref{lemmalattice} and \cite[Theorem 0.5.2]{dualmonoid},
the premonoid $P_c:=[1,c]$ (together with the natural partial product)
is a Garside premonoid.
We obtain a Garside monoid $\mathbf{M}(P_c)$.

By Proposition \ref{posetiso}, the restriction of $\pi$ is a bijection
from $S$ to $P_c$. Let $\phi$ be the inverse bijection.
Let $w,w'\in P_c$.
Assume that the product $w.w'$ is defined in $P_c$. Let $w''$ be the value
of this product. One has $w \preccurlyeq_R w''$, thus
$\phi(w) \preccurlyeq \phi(w'')$ (using again Proposition \ref{posetiso}),
and we may find $b'$ in $S$ such that $\phi(w) b' = \phi(w'')$ (Proposition
\ref{propiiiii}).
Claim: $b'= \phi(w')$. Indeed, $b'$ and $\phi(w')$ are two simple elements
whose image by $\pi$ is $w^{-1}w''$; one concludes using Corollary
\ref{corosimples}.
This proves that $\phi$ induces a premonoid morphism $P_c \rightarrow S$
(where $S$ is equipped with the restriction of the
monoid structure), thus induces a monoid
morphism $\mathbf{M}(P_c) \rightarrow M$ and 
a group morphism $\Phi:\mathbf{G}(P_c) \rightarrow B$.

Let use prove that $\Phi$ is an isomorphism.
Choose a basepoint $y\in Y^{\gen}$.
Let $\gamma_1,\dots,\gamma_n$ be generators of $\pi_1(Y-\CK,y)$.
Let $(s_1,\dots,s_n)$ be the label of $y$.
Let us reinterpret the presentation from
Theorem \ref{theovankampen}  in terms of Hurwitz action.
Since $LL:  Y-\CK \rightarrow E_n^{\gen}$ is a covering,
$\pi_1(Y-\CK,y)$ may be identified with a subgroup $H\subseteq B_n$.
The generators $f_1,\dots,f_n$ in Theorem \ref{theovankampen} may
be choosen to be $s_1,\dots,s_n$, and the monodromy automorphism
$\phi_1,\dots,\phi_m$ are obtained by Hurwitz action on $s_1,\dots,s_n$.
Let $h\in H$. Let $(s'_1,\dots,s'_n):= h(s_1,\dots,s_n)$ (by this,
we mean the Hurwitz action of $h$ on the free group generated by
$s_1,\dots,s_n$; the $s_i'$'s are words in the $s_i$'s).
Call \emph{Van Kampen relations} the relations of the
type $s'_i = s_i$, for any $i\in \{1,\dots,n\}$, $h\in H$, and $s'_i$
obtained as above.
We have
$$B\simeq \left< s_1,\dots,s_n \left|
\text{ Van Kampen relations} \right. \right>.$$
The map $\pi$ induces a bijection from 
 $A:=\pi^{-1}(R_c) \subseteq S$ to $R_c$.
Let $r_1,\dots,r_n$ be the images of $s_1,\dots,s_n$.
By transitivity of Hurwitz action on $\Red_R(c)$, the group
$\mathbf{G}(P_c)$ is generated by $r_1,\dots,r_n$, the remaining
generators in the presentation of Lemma \ref{lemmadualbraidrelations}
appearing as conjugates of $r_1,\dots,r_n$ (by successive use of dual
braid relations).
Our generating sets are compatible and the morphism
$$\Phi: 
\mathbf{G}(P_c) \simeq \left< R_c \left|
\text{dual braid relations} \right. \right> \rightarrow B \simeq
\left< s_1,\dots,s_n \left|
\text{ Van Kampen relations} \right. \right>$$
is defined by $r_i\mapsto s_i$.
Add to the presentation of $B$
 formal generators indexed by $$A -\{s_1,\dots,s_n\}$$
as well as the dual braid relations $\pi(r)\pi(r')=\pi(r')\pi(r'')$.
Since the relations already hold in $\mathbf{G}(P_c)$, they hold in $B$,
and we obtain a new presentation
$$B\simeq \left<A  \left|
\text{ Van Kampen relations on } \{s_1,\dots,s_n\}, \text{ dual braid
relations on } A \right. \right>.$$
To conclude that $\Phi$ is an isomorphism, it is enough to observe
that the dual braid relations encode the full Hurwitz action of
$B_n$ of $(s_1,\dots,s_n)$, while the Van Kampen relations encode the
action of $H\subseteq B_n$: thus Van Kampen relations are consequences
of dual braid relations, and $\mathbf{G}(P_c)$ and $B$ are given by
equivalent presentations.

Since $\Phi$ is an isomorphim and  $\mathbf{M}(P_c)$ naturally
embeds in $\mathbf{G}(P_c)$ (this is a crucial property of Garside
monoids, \cite{dehgar}), $\mathbf{M}(P_c)$ is isomorphic to its
image $M$ in $B$. The rest of theorem is clear.
\end{proof}

As mentioned earlier, one may view $B$ as a ``reflection group'', generated
by the set $\mathcal{R}$ of all braid reflections.
An element of $M$ is in $\mathcal{R}$ if and only if it has length $1$ for
the natural length function, or equivalently if it is an \emph{atom}
(\ie, an element which has no strict divisor in $M$ except the unity).
By Proposition \ref{propcardhurwitz} and Theorem \ref{theolooijenga},
there is a bijection between $\Red_R(c)$ and the image of $lbl:Y^{\gen}
\rightarrow \mathcal{R}^n$. This image is clearly contained
in $\Red_{\mathcal{R}}(\delta)$. The conjecture below claims an analog
in $B$ of the transitivity of the Hurwitz action on $\Red_R(c)$
 ($\delta$ is the natural substitute for a Coxeter
element in $B$).
It implies that any
element in  $\Red_{\mathcal{R}}(\delta)$ is the label of some $y\in Y^{\gen}$.

\begin{conj}
\label{braidreflections}
The Hurwitz action of $B_n$ on $\Red_{\mathcal{R}}(\delta)$ is transitive.
\end{conj}

The analog of this holds in free groups.

\section{Chains of simples}

Here again, $W$ is an irreducible well-generated complex $2$-reflection group,
and the notations from the previous section are still in use.
For $k = 0,\dots, n$, we denote by $\CC_k$ the set of (strict)
 chains in $S-\{1\}$ of
cardinal $k$, \ie, the set of $k$-tuples $(c_1,\dots,c_k)$ in $S^k$ such that
$$1 \prec c_1 \prec \dots \prec c_k,$$
or equivalently the set of $k$-tuples $(c_1,\dots,c_k)$ in $M^k$ such that
$$1 \prec c_1 \prec \dots \prec c_k \preccurlyeq \delta.$$
It is convenient to write $\{1 \prec c_1 \prec \dots \prec c_k \}$ instead
of $(c_1,\dots,c_k)$.

We set $\CC:=\bigsqcup_{k=0}^n \CC_k$.

Let $C:=\{1 \prec c_1 \prec \dots \prec c_k\} \in \CC$.
We say that $y\in Y$ \emph{represents}
$C$ if there exists $x\in U_y$ and real numbers $L_1,\dots,L_k$ such that
$0< L_1 < \dots < L_k$ and, for $i=1,\dots,k$, $(y,x,L_i)$ is
a tunnel representing $C_i$.

Example: if $LL(y)$ is as in the illustration below and
$(s_1,\dots,s_5) = lbl(y)$, $y$ represents
 $1\prec s_1\prec s_1s_3 \prec
s_1s_3 s_4\prec s_1s_3s_4s_5 $,
 $1\prec s_2 \prec s_2 s_4 \prec s_2s_4s_5$,
$1\prec s_2 \prec s_2s_5$
and their subchains,
but does not represent $1\prec s_2 \prec s_2s_3$
nor $1\prec s_3 \prec s_3s_5 $.
$$
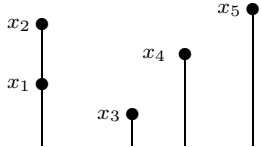
\xy
(-5,-10)="11",
(-5,-2)="1", (-5,6)="2", (7,-6)="3", (7,-10)="33",
(14,2)="4", (14,-10)="44",
(23,8)="5",(23,-10)="55",
"1"*{\bullet},"2"*{\bullet},"3"*{\bullet},"4"*{\bullet},
"5"*{\bullet},(20,8)*{_{x_5}}, "55";"5" **@{-},
(-8,-2)*{_{x_1}},(-8,6)*{_{x_2}},(4,-6)*{_{x_3}},(10,2)*{_{x_4}},
"11";"1" **@{-},
"2";"1" **@{-},
"33";"3" **@{-},
"44";"4" **@{-},
\endxy $$

\begin{defi}
For all $C\in \CC$, we set $Y_C := \{y\in Y | y \text{ represents }
C\}$.
\end{defi}

To illustrate this notion, we observe that implication (ii) $\Rightarrow$
(i) from Proposition \ref{propiiiii} expresses that for all $C\in \CC_2$,
$Y_C$ is non-empty. Based on the results from the
previous sections, this easily generalises to:

\begin{lemma}
For all $C\in \CC$, the space $Y_C$ is non-empty.
\end{lemma}

The goal of this section is to prove:

\begin{prop}
\label{YCcontractible}
For all $C\in \CC$, the space
$Y_C$ is contractible.
\end{prop}

This technical result will be used in Section \ref{sectionuniversalcover},
when studying the nerve of an open covering of the universal cover of
$W\cq V^{\reg}$: we will need to prove that certain non-empty intersections of
open sets are contractible, and these intersections will appear as
fiber bundles over some $Y_C$, with contractible fibers.

The proposition
is not very deep nor difficult but
somehow \emph{unconvenient} to prove since the retraction will be described
via $LL$, through ramification points. The following particular cases
are easier to obtain:
\begin{itemize}
\item If $C$ is the
chain $1 \prec \delta$, then $Y_C=Y \simeq
\BC^{n-1}$.
\item
More significatively, let $W$ be a complex reflection
group of type $A_2$. Up to renormalisation,
the discriminant is $X_2^2+X_1^3$. Identify $Y$ with $\BC$.
For all $y\in \BC$, $LL(y)=\{\pm (-y)^{3/2}\}$.
In particular, $LL(1) = \{\pm \sqrt{-1}\}$.
Let $s$ be the simple element represented by the tunnel
with $y=1$, $x=-1$ and $L=2$.
Let $C$ be the chain $1\prec s$.
Then $Y_C$ is the open cone consisting of non-zero elements
of $\BC$ with argument in the open interval $(-2\pi/3, 2\pi/3)$.
\item Assume that $C\in \CC_n$.
All points in $Y_C$ are
generic.
Consider the map $Y^{\gen} \rightarrow \CC_n\times E_n$
sending $y$ to the pair $(\{1\prec s_1
 \prec s_1s_2 \prec \dots \prec s_1s_2\dots s_n=\delta\},LL(y))$,
where $(s_1,\dots,s_n)=lbl(y)$. This map is an homeomorphism.
The $(Y_C)_{C\in \CC_n}$ are the connected components of $Y^{\gen}$. Each
of these components is homeomorphic to $E_n$, which is contractible.
These $(Y_C)_{C\in \CC_n}$ are some analogs of chambers.
\end{itemize}

In the following proposition, if $A\subseteq LL(y)$ is a submultiset,
the \emph{deep label} of $A$ is the sequence $(t_1,\dots,t_p)$
of labels (with respect to $y$) of points in the support of
 $A$ which are deep
with respect to $A$ (since these points may not be deep in $LL(y)$,
the deep label of $A$ is not necessarily a subsequence of the deep
label of $y$).

\begin{prop}
\label{inequalities}
Let $C= \{1\prec c_1 \prec c_2 \prec \dots \prec c_m \}$ be a chain
in $\CC_m$. Let $y\in Y$, let $(x_1,\dots,x_k)$ be the ordered support
of $LL(y)$ and $(s_1,\dots,s_k)$ be the label of $y$.
The following assertions are equivalent:
\begin{itemize}
\item[(i)] $y\in Y_C$.
\item[(ii)] There exists a partition of $LL(y)$ into
$m+1$ submultisets $A_0,\dots,A_m$ such that:
\begin{itemize}
\item[(a)] For all $i,j\in \{1,\dots,m\}$ with $i<j$, for all
$x\in A_i$ and all $x'\in A_j$, we have $\re(x) < \re(x')$.
\item[(b)] For all $x\in A_0$, one has
$$\re(x) < \min_{x'\in A_1\cup \dots \cup A_m} \re(x')$$ or
$$\im(x) < \min_{x'\in A_1\cup \dots \cup A_m} \im(x')$$ or
$$\re(x) > \max_{x'\in A_1\cup \dots \cup A_m} \re(x').$$
\item[(c)] For all $i\in \{1,\dots,m\}$, the product of
the deep label of $A_i$ is $c_i$.
\end{itemize}
\end{itemize}
Moreover, in situation (ii), the partition $LL(y) = A_0\sqcup\dots\sqcup A_m$
is uniquely determined by $y$ and $C$.
\end{prop}

The picture below illustrates the proposition for particular $y$ and $C$.
Here $lbl(y)=(s_1,\dots,s_5)$ and the considered chain is
$$C= \{1 \prec s_2 \prec s_2 s_4\}.$$
We have chosen tunnels $T_1=(y,x,L_1)$ and $T_2=(y,x,L_2)$, with $L_1< L_2$,
such that $b_{T_1}=s_2$ and $b_{T_2}=s_2s_4$.
The dotted lines represent these tunnels, as well as the 
vertical half-lines above $x$, $x+L_1$ and $x+L_2$. They partition
the complex line into three connected components; the partition
$A_0\sqcup A_1 \sqcup A_2$ is the associated partition of $LL(y)$.
It is clear the possibility of drawing the tunnels is subject precisely
to the conditions on $A_0\sqcup A_1 \sqcup A_2$ expressed in the
proposition.
$$
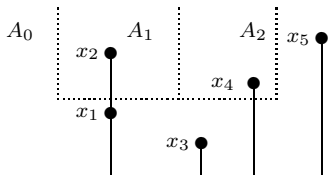
\xy
(-5,-10)="11",
(-5,-2)="1", (-5,6)="2", (7,-6)="3", (7,-10)="33",
(14,2)="4", (14,-10)="44",
(23,8)="5",(23,-10)="55",
(-17,9)*{_{A_0}},(-1,9)*{_{A_1}},(14,9)*{_{A_2}},
"1"*{\bullet},"2"*{\bullet},"3"*{\bullet},"4"*{\bullet},
"5"*{\bullet},(20,8)*{_{x_5}}, "55";"5" **@{-},
(-8,-2)*{_{x_1}},(-8,6)*{_{x_2}},(4,-6)*{_{x_3}},(10,2)*{_{x_4}},
"11";"1" **@{-},
"2";"1" **@{-},
"33";"3" **@{-},
"44";"4" **@{-},
(-12,0);(17,0) **@{.},
(-12,0);(-12,12) **@{.},
(4,0);(4,12) **@{.},
(17,0);(17,12) **@{.}
\endxy $$

\begin{proof}
(i) $\Rightarrow$ (ii):
$(y,z,L_1),\dots,(y,z,L_m)$ be tunnels representing the successive
non-trivial terms of $C=\{1\prec c_1 \prec \dots \prec c_m\}$.
 Set $z_0:=z$ and, for $i=1,\dots,m$, $z_i:=z+L_i$.
We have $$\re(z_0) < \re(z_1) < \dots < \re(z_m)$$ and
$$\im(z_0) = \im(z_1) = \dots = \im(z_m).$$
For $i\in\{1,\dots,m\}$, let $A_i$ be
the submultiset of $LL(y)$ consisting of points $x$ such that
$\re(z_{i-1}) < \re(x) < \re(z_i)$ and 
$\im(x) > \im(z_0)$. Let $A_0$ be the complement in $LL(y)$
of $A_1\cup \dots \cup A_m$. One easily checks (a), (b) and (c).

(ii) $\Rightarrow$ (i): Conversely, assume we are given a partition
$A_0\sqcup A_1\sqcup \dots \sqcup A_m$
satisfying conditions (a) and (b).
One may recover tunnels  $(y,z,L_1),\dots,(y,z,L_m)$ such that
the above construction yields the partition
$A_0\sqcup A_1\sqcup \dots \sqcup A_m$.
Condition (c) then implies that the tunnels
represent the elements of $C$, thus that $y\in Y_C$.

Unicity of the partition: this follows from condition (c) and
Lemma \ref{lemmaunique}.
\end{proof}

\begin{lemma}
\label{YC0}
Let $C= \{1\prec c_1 \prec c_2 \prec \dots \prec c_m \} \in \CC$.
Define $Y_C^0$ as the subspace of $Y_C$ consisting of points $y$ whose
associated partition $A_0,\dots,A_m$ (from Proposition \ref{inequalities}
(ii)) satisfies the following conditions:
\begin{itemize}
\item for $i=0,\dots,m$, the support of $A_i$ is a singleton $\{a_i\}$ and,
\item $\re(a_0) = \min_{i=1,\dots,m} \re(a_i) - 1$ and
$\im(a_0) = \min_{i=1,\dots,m} \im(a_i) - 1$.
\end{itemize}
Then $Y_C^0$ is contractible.
\end{lemma}

\begin{proof}
Let $y \in Y_C^0$. The support $(x_0,\dots,x_m)$
of $LL(y)$ is generic, thus the label $(s_0,\dots,s_m)$ of $y$ coincides
with the deep label, and we have $s_0\dots s_m=\delta$ (Corollary
\ref{corodeeplabel}).
By Proposition \ref{inequalities}, Condition (ii)(c),
we have, for $i=1,\dots,m$, $s_i=c_i$.
Thus $s_0= \delta (c_1\dots c_m)^{-1}$.
We have proved that the label of any $y\in Y_C^0$ must
be $(\delta (c_1\dots c_m)^{-1},c_1,\dots,c_m)$.
A consequence of Theorem \ref{theolabel} is that
 the map $Y_C^0\rightarrow E_{m}^{\gen}$ sending
$y$ to $(x_1,\dots,a_m)$ is an homeomorphism.
One concludes with Lemma \ref{lemmaEngen}.
\end{proof}

Regular coverings satisfy unique homotopy lifting property.
The map $LL: Y \rightarrow E_n$ is not regular on $\CK$ but
is \emph{stratified regular}, where $E_n$ is endowed with the
natural stratification: at each $y\in Y$, the Jacobian of $LL$
is surjective onto the tangent space of the stratum of $LL(y)$.
The following lemma says that paths
 of non-decreasing ramification admits unique lifts. It is a general
property of {stratifed regular} maps, though we only need
this particular example.

\begin{lemma}
\label{homotopylifting}
Let $\gamma:[0,1] \rightarrow E_n$ be a path.
For all $t\in [0,1]$, denote by $S_t$ the stratum of the natural
stratification containing $\gamma(t)$. Assume that $\forall t,t'\in[0,1],
t \leq t' \Rightarrow S_t \leq S_{t'}$. Let $y_0\in LL^{-1}(\gamma(0))$.
There exists a unique path $\tilde{\gamma}:[0,1]\rightarrow Y$ such that
$\gamma= LL\circ \tilde{\gamma}$ and $\tilde{\gamma}(0) = y_0$.
\end{lemma}

\begin{proof}
Left to the reader.
\end{proof}

We may now proceed to the proof of Proposition \ref{YCcontractible}.
Let $C=\{1\prec c_1 \prec \dots \prec c_m\} \in \CC$.
Let $y\in Y_C$. Let $(x_1,\dots,x_k)$ be the ordered support
of $LL(y)$, and $(n_1,\dots,n_k)$ be the multiplicities.
Let $A_0,\dots,A_m$ be the partition of $LL(y)$ described in 
Proposition \ref{inequalities} (ii).

The picture below gives an idea of the retraction of $Y_C$ onto $Y_C^0$ 
that will be explicitly constructed. It illustrates the motion of
a given point $y\in Y_C$; the black dots indicate the support of
$LL(y)$ and the arrows how this support moves during the retraction.
$$
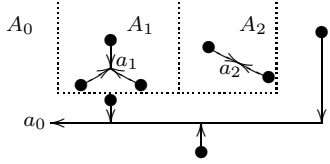
\xy
(-5,-4)="11",(-13,-4)="0",(-15,-4)*{_{a_0}},
(-5,-1)="1", (-5,7)="2", (-1,1)="22", (-9,1)="222",(-5,3.3)="20",
(-2.7,4)*{_{a_1}},
(7,-8)="3", (7,-4)="33",
(16,2)="4", (8,6)="44",(12,4)="40", (11,3)*{_{a_2}},
(23,8)="5",(23,-4)="55",
(-17,9)*{_{A_0}},(-1,9)*{_{A_1}},(14,9)*{_{A_2}},
"1"*{\bullet},"2"*{\bullet},"3"*{\bullet},"4"*{\bullet},
"22"*{\bullet},"222"*{\bullet},"44"*{\bullet},
"5"*{\bullet}, 
"2";"20" **@{-}, *\dir{>},
"22";"20" **@{-}, *\dir{>},
"222";"20" **@{-}, *\dir{>},
"4";"40" **@{-}, *\dir{>},
"44";"40" **@{-}, *\dir{>},
"1";"11" **@{-}, *\dir{>},
"3";"33" **@{-}, *\dir{>},
"5";"55" **@{-}, *\dir{>},
"55";"0" **@{-}, *\dir{>},
(-12,0);(17,0) **@{.},
(-12,0);(-12,12) **@{.},
(4,0);(4,12) **@{.},
(17,0);(17,12) **@{.}
\endxy $$

For $i=1,\dots,m$, consider the multiset mass centre
$$a_i:= \frac{\sum_{x \in A_i} x}{|A_i|}$$
(in this expression, $A_i$ is viewed as a multiset: each $x_j$ in $A_i$ is taken
$n_j$ times, and $|A_i|$ is the multiset cardinal, \ie, the sum of the 
$n_j$ such that $x_j\in A_i$).

For each $t\in [0,1]$, let
$$\gamma_y(t):= A_0\cup \bigcup_{i=1}^m \{ (1-t) x + ta_i | x\in A_i\}$$
(here again, we consider the multiset union; in particular, the multicardinal
of $\gamma_y(t)$ is constant, equal to $n$ -- \ie, $\gamma_y(t)\in E_n$).
This defines a path in $E_n$.

The path $\gamma_y$ satisfies the hypothesis
 of Lemma \ref{homotopylifting} and uniquely lifts to a path
$\tilde{\gamma}_y$ in $Y$ such that $\tilde{\gamma}_y(0)=y$. An easy
consequence of Proposition \ref{inequalities} is that
$\tilde{\gamma}_y$ is actually drawn in $Y_C$.

Let $y':=\tilde{\gamma}_y(1)$.
Let $$R:=\min_{i=1,\dots,m} \re(a_i) = \re(a_1)$$ and
$$I:=\min_{i=1,\dots,m} \im(a_i).$$
For all $x \in A_0$, we have $\re(x) < R$ or $\im(x) < I$ or
$\re(x) > \re(a_m)$; denote by $x'$ the complex number with the same real part
as $x$ and imaginary part $I-1$; let $a_0:=R-1+\sqrt{-1}(I-1)$.
Consider the path $\beta_x: [0,1] \rightarrow \BC$ defined by:
$$\beta_x(t) := \left\{ \begin{matrix} (1-2t)x + 2t x' &
\text{if } t\leq 1/2 \\ 
(2-2t) x' + (2t-1) a_0 &  \text{if } t \geq 1/2. \end{matrix} \right.$$
For all $t\in[0,1]$, one has $\re(\beta_x(t)) < R$ or
$\im(\beta_x(t)) < I$ or
$\re(\beta_x(t)) > \re(a_m)$.

The path $\gamma_y':[0,1] \rightarrow E_n$ defined by 
$$\gamma_y'(t):= 
\{\underbrace{a_1,\dots,a_1}_{|A_1| \text{ times}}\} \cup \dots \cup
\{\underbrace{a_m,\dots,a_m}_{|A_m| \text{ times}}\} \cup
\bigcup_{x\in A_0} \beta_x(t) $$
satisfies the hypothesis of Lemma \ref{homotopylifting} and lifts to
a unique path $\tilde{\gamma}_y'$ in $Y$ such that $\tilde{\gamma}_y'(0) = y'$.
Once again, a direct application of Proposition \ref{inequalities}
ensures that $\tilde{\gamma}_y'$ is actually in $Y_C$.
The endpoint $y'' := \tilde{\gamma}_y'(1)$ lies in the subspace
$Y_C^0$ of Corollary \ref{YC0}.

The map
\begin{eqnarray*}
\varphi:Y_C\times [0,1] & \longrightarrow & Y_C \\
(y,t) & \longmapsto & (\tilde{\gamma}_y' \circ \tilde{\gamma}_y)(t) 
\end{eqnarray*} is a retraction of $Y_C$ onto its
contractible subspace $Y_C^0$.
Thus $Y_C$ is contractible.

\section{The universal cover of $W \cq V^{\reg}$}
\label{sectionuniversalcover}

The main result of this section is:

\begin{theo}
\label{theouniversalcover}
The universal cover of $W \cq V^{\reg}$ is homotopy equivalent
to $\Flag(B,S)$.
\end{theo}

Combined with Corollary \ref{coroflag}, this proves our main result
Theorem \ref{theokapiun} in the irreducible case. The reducible case
follows.

To prove this theorem, we construct an open covering $\widehat{\CU} = 
(\widehat{\CU}_b)_{b\in B}$, such that the nerve
of $\CU$  is $\Flag(B,S)$ (Proposition \ref{nerve}) and 
intersections
of elements of $\widehat{\CU}$ are either empty or contractible (Proposition
\ref{intersections}). 
Under these assumptions, a standard theorem from algebraic topology
(\cite[4G.3]{hatcher}) gives the desired result.

As explained in the topological preliminaries, our ``basepoint'' $\CU$
provides us with a model denoted
$(\widehat{W \cq V^{\reg}})_{\CU}$ or simply $\widehat{W \cq V^{\reg}}$
for the universal cover of $W\cq V^{\reg}$.
Recall that, to any semitunnel $T$, one associates a path $\gamma_T$ whose
source is in $\CU$, thus a point in $\widetilde{W\cq V^{\reg}}$.
Moreover, we have a left action of $B=\pi_1(W \cq V^{\reg}, \CU)$
on $\widetilde{W \cq V^{\reg}}$.
With these conventions, our open covering is very easy to define:

\begin{defi}
The set $\widehat{\CU_1}$ is the subset of $\widetilde{W \cq V^{\reg}}$
of elements represented by semitunnels.
For all $b\in B$, we set $\widehat{\CU}_b:= b \widehat{\CU}_1.$
\end{defi}

It is clear that two semitunnels represent the
same point in $\widetilde{W \cq V^{\reg}}$ if and only if they are
equivalent in the following sense:

\begin{defi}
Two semitunnels $T=(y,x,L)$ and $T'=(y',x',L')$ are \emph{equivalent}
if and only if $y=y'$, $x+L=x'+L'$ and the affine segment $[(x,y),(x',y)]$ is
included in $\CU$.
\end{defi}

Let $T=(y,x,L)$ be a semitunnel.
The point of $\widehat{\CU}_1$ determined by $T$, or in other
words the equivalence class of $T$, is uniquely determined by $y$, $x+L$ and
$$\lambda(T):=\inf \left\{ \lambda'\in [0,L]
\left| [(x,y),(x+L-\lambda',y)] \subseteq \CU\right.
\right\}.$$
The number $\lambda$ is the infimum of the length of semitunnels in the
equivalence class of $T$. This infimum is not a minimum (unless
$T$ is included in $\CU$), since
$(y,x+L-\lambda(T),\lambda(T))$ is not in general a semitunnel.

We consider the following subsets of 
$\widehat{\CU}_1$:
\begin{itemize}
\item[($\CU_1$)] If 
$\lambda(T)=0 =\min \left\{ \lambda'\in [0,L] \left|
[x,x+L-\lambda'] \subseteq U_y\right.
\right\}$, then $T$ is a tunnel, equivalent to $(y,x+L,0)$. Elements
of $\widehat{\CU}_1$ represented by such tunnels of length $0$
form an open subset denoted by $\CU_1$. This subset is actually a sheet
over $\CU$ of the universal covering, corresponding to the trivial lift
of the ``basepoint''.
\item[($\overline{\CU}_1$)]
We denote by $\overline{\CU}_1$ the subset consisting of points represented
by tunnels with $\lambda(T)=0$ (but without requiring that the
 ``$\inf$'' is actually a ``$\min$'').
We obviously have $\CU_{1} \subseteq \overline{\CU}_1$,
and $\overline{\CU}_1$ is contained in the closure of $\CU_1$.
\end{itemize}
Similarly, we set for all $b\in B$
$$\CU_b:=b\CU_1, \quad \overline{\CU}_b:=b\overline{\CU}_1.$$
Also, for any $y\in Y$, we denote by 
$\widehat{\CU}_{b,y}$
(resp. $\CU_{b,y}$, resp. $\overline{\CU}_{b,y}$) the intersection of 
$\widehat{\CU}_{b}$
(resp. $\CU_{b}$, resp. $\overline{\CU}_{b}$)
with the fiber over $y$ of the composed map
$\widetilde{W\cq V^{\reg}} \rightarrow W \cq V^{\reg} \rightarrow Y$.

\begin{lemma}
The family $(\overline{\CU}_b)_{b\in B}$ is a partition of
$\widetilde{W\cq V^{\reg}}$.
\end{lemma}

\begin{proof}
It suffices to show that the projection $\overline{\CU}_1 \rightarrow
W \cq V^{\reg}$ is bijective.
Any point in $z\in W\cq V^{\reg}$ is the target of a unique equivalence
class of semitunnels $T$
with $\lambda(T)=0$; depending on whether $z\in \CU$ or not, the
associated point will be in $\CU_1$ or in $\overline{\CU}_1- \CU_1$.
\end{proof}

In particular, $(\widehat{\CU}_b)_{b\in B}$ is a covering of
$\widetilde{W\cq V^{\reg}}$.

\begin{lemma}
\label{lemmacon}
For all $b\in B$, $\widehat{\CU}_b$ is open and contractible.
\end{lemma}

\begin{proof}
It is enough to deal with $b=1$.
That $\widehat{\CU}_1$ is open is easy.
Let $\CT$ be the space of tunnels and $\sim$ the equivalence relation.
As a set, $\widehat{\CU}_1\simeq \CT/\sim$.
Consider the 
map
\begin{eqnarray*}
\phi:  \CT \times [0,1] & \longrightarrow  & \CT  \\
(T=(y,x,L),t) & \longmapsto & (y,x,L-t\lambda(T)) 
\end{eqnarray*}
If $T\sim T'$, then $\forall t, \phi(T,t)\sim \phi(T',t)$.
Thus $\phi$ induces a map $\overline{\phi}:\widehat{\CU}_1
\times [0,1] \rightarrow \widehat{\CU}_1$.
This map is continuous (this follows
from the fact that $\lambda$ induces a continuous function on
$\widehat{\CU}_1$) and $\overline{\phi} (T,t)=T$ if $T\in
\overline{\CU}_1$ or if $t=0$. We have proved that
$\widehat{\CU}_1$ retracts to $\overline{\CU}_1$.
We are left with having to prove that $\overline{\CU}_1$ is contractible.
Inside $\overline{\CU}_1$ lies $\CU_1$ which is contractible, since it
is a standard lift of the contractible space $\CU$.

We conclude by observing that
$\CU_1$ and $\overline{\CU}_1$ are homotopy equivalent. There
is probably a standard theorem from semialgebraic geometry applicable
here, but I was unable to
find a proper reference. Below is a ``bare-hand'' argument: it explains
how $\overline{\CU}_1$ may be ``locally retracted'' inside $\CU_1$
(constructing a global retraction seems difficult).

Both spaces have homotopy
type of $CW$-complexes and to prove homotopy equivalence it is enough
to prove that any continuous map $f:S^k \rightarrow \overline{\CU}_1$
(where $S^k$ is a sphere)
may be homotoped to a map $S^k \rightarrow {\CU}_1$.
Assume that there is a tunnel $T=(y,z,L)$ such that $(T/\sim) \in f(S^k)
\cap (\overline{\CU}_1 - \CU_1)$.
We have $L>0$.
For any $\varepsilon > 0$, let $B_\varepsilon(y)$ be the open
ball of radius $\varepsilon$ in $Y$ around $y$, let $B_\varepsilon(z)$
be the affine interval $(z-\sqrt{-1}\varepsilon,z+\sqrt{-1}\varepsilon)$.
For $\varepsilon$ small enough, there is a unique continuous function
$L_{\varepsilon}: B_\varepsilon(y)\times B_\varepsilon(z) \rightarrow \BR$
such that $L_\varepsilon(y,z)=L$ and for all $(y',z')\in
 B_\varepsilon(y)\times B_\varepsilon(z)$, $(y',z',L_{\varepsilon}
(y',z'))$ represents a point in $\overline{\CU}_1 - \CU_1$.
The ``half-ball''
$$H_\varepsilon := \{(y',z',L') \in B_\varepsilon(y)\times
 B_\varepsilon(z) \times \BR | 0\leq L' \leq L_{\varepsilon}(y',z') \}$$
is a neighbourhood of $T/\sim$ in $\overline{\CU_1}$.
Working inside this neighbourhood, one may homotope $f$ to $f'$ such
that $f'(S^k)
\cap (\overline{\CU}_1 - \CU_1) \subseteq f(S^k)
\cap (\overline{\CU}_1 - \CU_1)-\{T/\sim\}$.
Compacity of $f(S^k)
\cap (\overline{\CU}_1 - \CU_1)$ garantees that one can iterate this
process a finite number of times to get rid of all $f(S^k)
\cap (\overline{\CU}_1 - \CU_1)$.
\end{proof}

\begin{prop}
\label{nerve}
The nerve of $(\widehat{\CU}_b)_{b\in B}$ is
$\Flag(B,S)$.
\end{prop}

\begin{proof}
Let $b,b'\in B$ such that
$\widehat{\CU}_b\cap \widehat{\CU}_{b'}\neq \varnothing$.
Let $T=(y,x,L)$ and $T'=(y',x',L')$ be semitunnels, representing
points $z$ and $z'$ in $\widehat{\CU}_1$, such that
$bz=b'z'$. The image of $z$ (resp. $z'$) in $W\cq V^{\reg}$ is
$(x+L,y)$ (resp. $(x'+L',y')$). Thus $x+L=x'+L'$ and $y=y$.
Up to permuting $b$ and $b'$, we may assume that $L\geq L'$.
Since $x'=x+L-L'$ is in $U_y$, $T'':=(y,x,L-L')$ is a tunnel,
representing a simple element $b''$.
The tunnel $T$ is a concatenation of $T''$ and $T'$.
This implies that $z=b''z'$ and $b'z'=bz=bb''z'$.
By faithfullness of the $B$-action on the orbit of $z$,
we conclude that $bb''=b'$.

We have proved that the $1$-skeletons of the nerve and of
$\Flag(B,S)$ coincide. To conclude, it remains to
check that the nerve is a flag complex. Let $C\subseteq B$ be such
that for all $b,b'\in C$, either $b^{-1}b'$ or $b'^{-1}b$ is simple.
We have to prove that $\bigcup_{b\in C} \widehat{\CU}_b\neq
\varnothing.$

Let $c_0,\dots,c_m$ be the elements of $C$, numbered according
to the total ordering on $C$ induced by $\preccurlyeq$:
$$c_0 \prec c_1 \prec c_2 \prec \dots \prec c_m \preccurlyeq c_0\delta.$$
Up to left-dividing each term by $c_0$, we may assume that $c_0=1$.
Let $y\in Y_C$. We may find $x,L_1,\dots,L_m$ such
that $(y,x,L_i)$ represents $c_i$ for all $i$.
The point of the universal cover represented by $(y,x,L_m)$
belongs to $\bigcup_{b\in C} \widehat{\CU}_b$.
\end{proof}

\begin{prop}
\label{intersections}
Let $C$ be a subset of $B$ such that
 $\bigcap_{b\in C}\widehat{\CU}_b \neq \varnothing $. Then 
 $\bigcap_{b\in C}\widehat{\CU}_b$ is contractible.
\end{prop}

\begin{proof}
As in the previous proof, we write $C=\{c_0,c_1,\dots,c_m\}$
with $$c_0 \prec c_1 \prec c_2 \prec \dots \prec c_m \preccurlyeq c_0\delta$$
and assume without loss of generality that $c_0=1$.

The case $m=0$ is Lemma \ref{lemmacon}.

Assume that $m\geq 1$. 
Let $T=(y,z,L)$. The point represented by
$T$ lies in $\bigcap_{b\in C}\widehat{\CU}_b$ if and only
if there exists $L_1,\dots,L_m$ with $0< L_1 < \dots < L_m < L$ such that,
for all $i$, $T_i:=(y,z,L_i)$ is a tunnel representing $c_i$.
Given $y\in Y$, it is possible to find such $L_1,\dots,L_m$ if and only
if $y\in Y_C$. This justifies:
$$\bigcap_{b\in C}\widehat{\CU}_{b,y} \neq \varnothing \Leftrightarrow
y \in Y_C.$$

For a given $y\in Y_C$, let us study
the intersection $\bigcap_{b\in C}\widehat{\CU}_{b,y}$.
Let $(x_1,\dots,x_k)$ be the ordered support of $LL(y)$.
Let $A_0,\dots,A_m$ be the associated partition of $LL(y)$,
defined in Proposition \ref{inequalities}.
Let $$R_+(y,C):=\max\{ \re(z) | z\in A_m\},$$
$$R_-(y,C):=\min\{ \re(z) | z\in A_1\},$$
 $$I_+(y,C):=\min\{\im(z) | z\in A_1\cup\dots \cup A_m\},$$
and
 $$I_-(y,C)=\sup\{\im(z) | z \in A_0 \text{ and } R_- \leq \re(z) \leq R_+\}.$$
We have $I_-(y,C) < I_+$. It may happen that $I_-(y,C)=-\infty$. 

We illustrate this on a picture. In the example,
the support of $LL(y)$ is $x_1,\dots,x_6$
and $C=\{1\prec s_2 \prec s_2s_4\}$, where $(s_1,\dots,s_6)=lbl(y)$.
The support of $A_1$ is $x_2$ and the support of $A_2$ is $x_4$.
The remaining points are in $A_0$.
The lines $\im(z)=I_-(y,C)$ and $\im(z)=I_+(y,C)$ are represented by full lines.
A semitunnel $(y,z,L)$
representing a point in $\bigcap_{b\in C}\widehat{\CU}_{b,y}$
must cross the intervals $I_2$ and $I_4$ represented by dotted lines.
One must have $\re(z) < R_-(y,C)$ and $I_-(y,C) < \im(z) < I_+(y,C)$; the final
point $z+L$ must satisfy $\re(z+L) > R_+(y,C)$. This final point
 may be any complex number
$z'$ in the rectangle $\re(z') > R_+(y,C)$ and $I_-(y,C) < \im(z')
< I_+(y,C)$, except
the points on the closed horizontal half-line to the right of $x_6$ (indicated
by a dashed line), which
cannot be reached.
$$
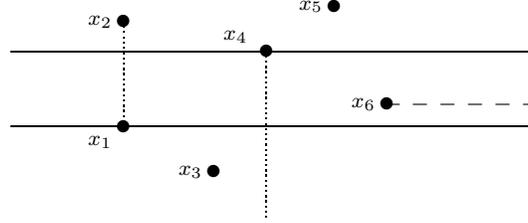
\xy
(-5,-10)="11",
(-5,-8)="1", (-5,6)="2", (7,-14)="3", (7,-10)="33",
(14,2)="4", (14,-20)="44",
(23,8)="5",(30,-5)="6",(50,-5)="66",
"1"*{\bullet},"2"*{\bullet},"3"*{\bullet},"4"*{\bullet},
"5"*{\bullet},(20,8)*{_{x_5}},
"6"*{\bullet},(27,-5)*{_{x_6}},
(-8,-10)*{_{x_1}},(-8,6)*{_{x_2}},(4,-14)*{_{x_3}},(10,4)*{_{x_4}},
"6";"66" **@{--},
(-20,-8);(50,-8) **@{-},
(-20,2);(50,2) **@{-},
"1";"2" **@{.},
"4";"44" **@{.}
\endxy $$

Generalising the example, one shows that $\bigcap_{b\in C}\widehat{\CU}_{b}$
may be identified with $$\bigcup_{y \in Y_C} E(y,C),$$
where $E(y,C)$ is the open rectangle of $\BC$ defined by
$$\re(z) > R_+(y,C), \quad I_-(y,C) < \im(z) < I_+(y,C),$$
from which have been removed the possible points of $A_0$ and the
horizontal half-lines to their rights.

Let $\overline{E}(y,C)$ be the rectangle
$$\re(z) \geq R_+(y,C), \quad I_-(y,C) < \im(z) < I_+(y,C)$$
from which have been removed the possible points of $A_0$ and the
horizontal half-lines to their rights. 

A homotopy argument similar to the one used in the proof of Lemma
\ref{lemmacon} (or possibly a nicer argument from semialgebraic geometry)
shows that $\bigcup_{y \in Y_C} E(y,C)$ and $\bigcup_{y \in Y_C} \overline{E}
(y,C)$ are homotopy equivalent.
The latter may be retracted to the union of open intervals
$\bigcup_{y\in Y_C} (I_-(y,C),I_+(y,C))$
(on each rectangle, the retraction is 
$$\overline{E}(y,C) \times [0,1] \rightarrow  \overline{E}(y,C),
(z,t) \mapsto R_+(y,C) + t(\re(z) -R_+(y,C)) 
+ \sqrt{-1} \im(z)).$$

The union $\bigcup_{y\in Y_C} (I_-(y,C),I_+(y,C))$ may be retracted
to $\bigcup_{y\in Y_C} [I_0(y,C),I_+(y,C))$, where
$I_0(y,C):= \max(\frac{I_-(y,C) +I_+(y,C)}{2},I_+(y,C)-1)$. The latter
space is a fibre bundle over $Y_C$, with fibers intervals.
Since the basespace $Y_C$ is contractible (Proposition \ref{YCcontractible}),
this fibre bundle is contractible. So is $\bigcap_{b\in C}\widehat{\CU}_b$.
\end{proof}

\section{Generalised non-crossing partitions}

Here again, $W$ is an irreducible well-generated complex $2$-reflection 
group.

When $W$ is of type $A_{n-1}$, the lattice $(S,\preccurlyeq)$ is isomorphic
to the lattice of non-crossing partitions of a regular $n$-gon
(\cite{BDM}, \cite{brady}).
Following \cite{reiner}, \cite{dualmonoid} and \cite{BC}, we call \emph{generalised lattice of non-crossing partitions of type $W$}
the lattice $$(S,\preccurlyeq),$$
 and \emph{Catalan number of type $W$} 
the number
$$\Cat(W):= \prod_{i=1}^n \frac{d_i+h}{d_i}.$$

The operation sending $s\preccurlyeq t$ to $s^{-1}t$ is an analog
of Kreweras complement operation. The map $s\mapsto s^{-1}\delta$ is
an anti-automorphism of the lattice.

In the Coxeter case,
Chapoton discovered a general formula for the number
$Z_W(N)$ of weak chains $s_1\preccurlyeq s_2 \preccurlyeq \dots
\preccurlyeq s_N$ of length $N$ in $(S,\preccurlyeq)$ (see \cite{chapoton}).
This formula continues to hold, though we are only able to prove this
case-by-case
(see \cite{athanasiadisreiner} and \cite{chapoton} for the Coxeter types;
the $G(e,e,n)$ case was done in \cite{BC};
the remaining types are done by
computer).

\begin{prop}
We have, for all $N$, 
$$Z_{W}(N)=\prod_{i=1}^n \frac{d_i+ N h}{d_i}.$$
\end{prop}

\begin{coro}
We have $|S| = \Cat(W)$.
\end{coro}

Another interesting numerical invariant is the \emph{Poincar\'e polynomial}
$$\Poin(S) := \sum_{s\in S} t^{l(s)}.$$

The table below summarises the numerical data for the exceptional types
(real and non-real). The coefficient of $t$ in the Poincar\'e polynomial
is the cardinal of $R_c$. One observes that $$R=R_c \Leftrightarrow W
\text{ is real.}$$

In the Weyl group case, $\Poin(S)$ may be interpreted as the Poincar\'e
polynomial of the cohomology of a toric variety related to cluster algebras
(\cite{chapoton}).

\begin{tabular}{|c|c|c|c|c|c|}
\hline
$W$ & degrees & $|R|$ & $\Cat(W)$ & $\Poin(S)$ & $|\Red_T(c)|$ \\
\hline
$G_{23}$ ($H_3$) & $2,6,10$  & $15$ & $32$ & $1 + 15t + 15t^2 + t^3$ & $50$ \\
\hline
$G_{24}$ & $4,6,14$  & $21$ & $30$ & $1+14t+14t^2+t^3$ & $49$ \\
\hline
$G_{27}$ & $6,12,30$  & $45$ & $42$ & $1+20t+20t^2+t^3$ & $75$ \\
\hline
$G_{28}$ ($F_4$) & $2,6,8,12$ &  $24$ & $105$ & 
$1+24t+55t^2+24t^3+t^4$ & $432$ \\ 
\hline
$G_{29}$ & $4,8,12,20$ &  $40$ & $112$ & $1+25t+60t^2+25t^3+t^4$
& $500$ \\
\hline
$G_{30}$ ($H_4$) & $2,12,20,30$  & $60$ & $280$ &
$1+60t+158t^2+60t^3+t^4$ & $1350$\\
\hline
$G_{33}$ & $4,6,10,12,18$ &  $45$ & $308$ &
$\begin{matrix}1+30t+123t^2 \\+123t^3+30t^4+t^5\end{matrix}$ & $4374$ \\
\hline
$G_{34}$ & $\begin{matrix} 6,12,18,\\ 24,30,42\end{matrix}$ &
$126$  &  $1584$ &
$\begin{matrix}1+56t+385t^2+ 700t^3 \\
+385t^4+56t^5+t^6\end{matrix}$ & $100842$ \\
\hline
$G_{35}$ ($E_6$) & $\begin{matrix}2,5,6, \\8,9,12\end{matrix}$ & 
 $36$  &  $833$ & 
$\begin{matrix}1+36t+204t^2+351t^3\\ +204t^4+36t^5+t^6\end{matrix}$ & 
$41472$ \\
\hline
$G_{36}$($E_7$) & $\begin{matrix} 2,6,8,10,\\ 12,14,18\end{matrix}$ & $63$ &
$4160$ &
$\begin{matrix}1+ 63t+ 546t^2+ 1470t^3\\
+ 1470t^4+ 546t^5+ 63t^6 + t^7 \end{matrix}$ & 
$1062882$ \\
\hline
$G_{37}$ ($E_8$)& $\begin{matrix}2,8,12,14, \\18,20,24,30 \end{matrix}$ &
$120$
& $25080$ &
$\begin{matrix} 1+120t+1540t^2 \\
+6120t^3+9518t^4+6120t^5 \\ +1540t^6+120t^7+t^8\end{matrix}$ &
$37968750$ \\
\hline
\end{tabular}

\section{On badly-generated groups}

Our theory works with well-generated groups, but we conjecture
that the dual braid monoid is a particular
case of a more general construction, giving good monoids for
braid groups of badly-generated groups.
If $d=d_{i_0}$ is a
regular degree of an arbitrary irreducible complex reflection
group, let $L$ be a generic line of direction $X_{i_0}$ in $W\cq V$
(as in \cite[section 3]{zariski}). It intersects the discriminant 
variety $\CH$ in $(N+N^*)/d$ points. 
Let $C$ be an Euclidean circle in $L$ surrounding this intersection.
It represents a ``Coxeter element'' in the free group
$\pi_1(L-L\cap \CH)$, as in \cite{free}, and we have a natural
``dual braid monoid'' for the free group, generated by
non-crossing loops drawn inside $C$. We conjecture that the direct
image of this monoid in $B(W)$ is a \emph{quasi-Garside monoid}, \ie,
it satisfies all axioms of Garside monoids except the finiteness of the
set of simples (the finiteness should occur only if $W$ is well-generated
and $d=h$ -- it is not difficult to see that, in this case, this monoid
coincides with the one studied in the present article).
The main difficulty is that this prevents the use of a counting
argument such as the one used in this article.

\end{document}